\documentclass{amsart}
\usepackage{color}
\usepackage{amsmath}
\usepackage{kotex}
\usepackage{amssymb}
\usepackage{amsthm}
\usepackage{amscd}
\usepackage{bm}
\usepackage{eucal} 
\usepackage{enumerate}
\usepackage{tikz}
\usepackage{tikz-cd}
    \usetikzlibrary{positioning, fit, calc}
\usepackage{shuffle}
\usepackage[utf8]{inputenc}
\usepackage[all]{xy}
\allowdisplaybreaks

\usepackage{hyperref}
\hypersetup{
  colorlinks   = true,          
  urlcolor     = blue,          
  linkcolor    = blue,          
  citecolor   = red             
}

\theoremstyle{plain}
\newtheorem{theorem}{Theorem}[section]
\newtheorem{conjecture}[theorem]{Conjecture}
\newtheorem{lemma}[theorem]{Lemma}
\newtheorem{proposition}[theorem]{Proposition}

\theoremstyle{definition}
\newtheorem{definition}[theorem]{Definition}

\newtheorem{notation}[theorem]{Notation}

\numberwithin{equation}{section}

\newcommand{\ppar}{${}$\par}

\newcommand\fantome[1]{}

\def\bN{\mathbb N}

\def\fa{\mathfrak{a}}
\def\fb{\mathfrak{b}} 
\def\fc{\mathfrak{c}} 
\def\fu{\mathfrak{u}} 
\def\fv{\mathfrak{v}} 

\def\Fq{\mathbb F_q}

\newcommand{\sha}{\shuffle}

\DeclareMathOperator{\Id}{Id}

\newcommand{\F}{\mathbb{F}}

\newcommand{\fs}{\mathfrak{s}}
\newcommand{\ft}{\mathfrak{t}}  

\newcommand{\fve}{\bm{\varepsilon}}

\newcommand{\N}{\ensuremath \mathbb{N}}

\DeclareMathOperator{\depth}{depth}

\author[B.-H. Im]{Bo-Hae Im}
\address{
Dept. of Mathematical Sciences, KAIST,
291 Daehak-ro, Yuseong-gu,
Daejeon 34141, South Korea
}
\email{bhim@kaist.ac.kr}

\author[H. Kim]{Hojin Kim}
\address{
Dept. of Mathematical Sciences, KAIST,
291 Daehak-ro, Yuseong-gu,
Daejeon 34141, South Korea
}
\email{hojinkim@kaist.ac.kr}

\author[K. N. Le]{Khac Nhuan Le}
\address{
Normandie Université,
Université de Caen Normandie - CNRS,
Laboratoire de Mathématiques Nicolas Oresme (LMNO), UMR 6139,
14000 Caen, France.
}
\email{khac-nhuan.le@unicaen.fr}

\author[T. Ngo Dac]{Tuan Ngo Dac}
\address{
Normandie Université,
Université de Caen Normandie - CNRS,
Laboratoire de Mathématiques Nicolas Oresme (LMNO), UMR 6139,
14000 Caen, France.
}
\email{tuan.ngodac@unicaen.fr}

\author[L. H. Pham]{Lan Huong Pham}
\address{
Institute of Mathematics, Vietnam Academy of Science and Technology, 18 Hoang Quoc Viet, 10307 Hanoi, Viet Nam
}
\email{plhuong@math.ac.vn}

\title[Hopf algebras and AMZV's]{Hopf algebras and alternating multiple zeta values in positive characteristic}



\date{\today}

\begin{document}
\begin{abstract}
In \cite{IKLNDP23} we presented a systematic study of algebra structures of multiple zeta values in positive characteristic introduced by Thakur as analogues of classical multiple zeta values of Euler. In this paper we construct algebra and Hopf algebra structures of alternating multiple zeta values introduced by Harada, extending our previous work. Our results could be considered as an analogue of those of Hoffman \cite{Hof00} and Racinet \cite{Rac02} in the classical setting. The proof is based on two new ingredients: the first one is a direct and explicit construction of the shuffle Hopf algebra structure, and the second one is the notion of horizontal maps.
\end{abstract}

\maketitle
\tableofcontents

\section{Introduction}

\subsection{Classical setting} ${}$\par

Let $\mathbb N = \{1, 2, \dots\}$ be the set of positive integers. The multiple zeta values (MZV's for short) are infinite sum defined by
\[ \zeta(n_1, \dots, n_r) = \sum_{0 < k_1 < \dots < k_r} \frac{1}{k_1^{n_1} k_2^{n_2} \dots k_r^{n_r}},\]
for $n_i \in \mathbb{N}$ with $n_r \ge 2$ for convergence. Here $r$ and $w:= n_1 + \dots + n_r$ are said to be the depth and weight of the presentation $\zeta(n_1, \dots, n_r)$, respectively. In the 18th century, Euler \cite{Eul76} studied MZV's for depth 1 and 2 cases. After two centuries, Zagier~\cite{Zag94} studied the general MZV's, initiating extensive works connected to various fields in mathematics and physics, including number theory, K-theory, knot theory, and higher energy physics (see for example \cite{Bro12,DG05,Dri90,Hof97,Kon99,IKZ06,LM96,Ter02}).

There exist more general objects called colored multiple zeta values. Letting $N \in \mathbb N$ be a positive integer and $\Gamma_N$ be the group of $N$-th roots of unity, we define colored multiple zeta values by
	\[ \zeta \begin{pmatrix}
	\epsilon_1 & \dots & \epsilon_r \\
	n_1 & \dots & n_r
	\end{pmatrix}=\sum_{0<k_1<\dots<k_r} \frac{\epsilon_1^{k_1} \dots \epsilon_r^{k_r}}{k_1^{n_1} \dots k_r^{n_r}} \]
for $\epsilon_i \in \Gamma_N$ and $n_i \in \N$ with $(n_r,\epsilon_r) \neq (1,1)$ for convergence. These objects were first systematically studied by Arakawa-Kaneko, Deligne, Goncharov and Racinet (see \cite{AK99, AK04, Del10, DG05, Gon01, Gon05, Rac02}). Since then they have been also studied by a lot of mathematicians and physicists, including Broadhurst \cite{Broa96a, Broa96b}, Hoffman \cite{Hof19}, Kaneko-Tsumura \cite{KT13}, and Zhao \cite{Zha16} in various fields of mathematics and physics. For $N=1$ we recover MZV's. For $N=2$ the colored multiple zeta values are also known as alternating multiple zeta values (AMZV's for short) or Euler sums.

The main goal of the theory of colored multiple zeta values is to find all linear and algebraic relations among colored multiple zeta values over a suitable coefficient field. It led to several important conjectures formulated by Ihara-Kaneko-Zagier \cite{IKZ06}, Zagier \cite{Zag94} and Hoffman \cite{Hof97}. We will explain below Ihara-Kaneko-Zagier's conjecture and refer the reader to the excellent book \cite{BGF} for more details about Zagier and Hoffman's conjectures.

For MZV's (the case $N=1$) Ihara-Kaneko-Zagier \cite{IKZ06} introduced the double shuffle relations and regularization. Roughly speaking, they used Hoffman's algebra~$\mathfrak h$ (see \cite{Hof00}) and subalgebras $\mathfrak h^0 \subset \mathfrak h^1 \subset \mathfrak h$ which are related to the MZV's, and presented two different algebra structures $(\mathfrak h^1, *)$ and $(\mathfrak h, \sha)$. These structures can be extended to Hopf algebra structures (see for example the work of Hoffman \cite{Hof00}). By regularization one gets the double shuffle relations by `comparing' two products in the level of MZV's. Double shuffle relations and regularization produce linear relations among MZV's, and Ihara-Kaneko-Zagier conjectured that these relations are enough to understand all of the linear relations among MZV's (see \cite[Conjecture 1]{IKZ06}).

For colored multiple zeta values the previous works of Ihara-Kanako-Zagier and Hoffman were extended by Racinet in \cite{Rac02}. However, we warn the reader that as mentioned in {\it loc. cit.}, for general $N$ the double shuffle relations and regularization fail to exhaust all the linear relations among colored multiple zeta values.

\subsection{Positive characteristic setting}
\subsubsection{MZV's and AMZV's in positive characteristic} ${}$\par

Carlitz \cite{Car35} and Thakur \cite{Tha04, Tha17} defined the zeta and multiple zeta values in positive characteristic via the function fields analogy (see \cite{Iwa69,MW83,Wei39}). Let $q$ be a power of prime $p$ and $\Fq$ be a finite field of order $q$. We use the following notations:
\begin{align*}
    A = \Fq[\theta]\quad &\text{the ring of polynomials in the variable $\theta$ over $\F_q$,}\\
    A_+ \quad&\text{the set of monic polynomials in $A$,}\\
    K = \Fq(\theta)\quad &\text{the field of fractions of $A$ (with the rational point $\infty$),}\\
    K_\infty = \Fq(\!(1/\theta)\!) \quad &\text{the completion of $K$ at $\infty$,}\\
    v_\infty\quad&\text{the discrete valuation on $K$ associated to the place $\infty$},\\
    &\text{normalized as $v_\infty(\theta) = -1$},\\
    |\cdot|_\infty = q^{-v_\infty(\cdot)} \quad&\text{the corresponding absolute value on $K$.}
\end{align*}

Carlitz \cite{Car35} defined the Carlitz zeta values  $\zeta_A(n)$ for $n \in \N$ as
\[\zeta_A(n):= \sum_{a \in A_+} \frac{1}{a^n} \in K_\infty.\] Note that this definition is an analogue of the classical zeta values $\zeta(n) = \sum_{k\in\N} \frac{1}{k^n}$. Extending this, Thakur \cite{Tha04} defined the multiple zeta value in positive characteristic for tuples $\fs=(s_1, \dots, s_r) \in \N^r$ as
\[ \zeta_A(\fs) = \sum\frac{1}{a_1^{s_1} \dots a_r^{s_r}}.\]
Here the sum is over all tuples of monic polynomials with decreasing degree, i.e., $(a_1, \dots, a_r) \in A_+^r$ with $\deg(a_1)> \dots > \deg(a_r)$. One can see the analogy between the classical MZV's and MZV's in positive characteristic. We define the depth and weight of $\zeta_A(\fs)$ by $r$ and $w(\fs) := s_1 + \dots+ s_r$ respectively. Note that Carlitz zeta values are MZV's of depth one. Thakur showed some of their properties including non-vanishing \cite{Tha09a}, sum-shuffle relations \cite{Tha10}, etc. For further references on the MZV's in positive characteristic, see \cite{AT90,AT09,GP21,Gos96,LRT14,LRT21,Pel12,Tha04,Tha09,Tha10,Tha17,Tha20,Yu91}.

In \cite{Har21} Harada introduced the notion of alternating multiple zeta values in positive characteristic.
Let $\fs=(s_1,\dots,s_r) \in \N^r$ and $\fve=(\varepsilon_1,\dots,\varepsilon_r) \in (\Fq^*)^r$, we define
\begin{equation*}
    \zeta_A \begin{pmatrix}
 \fve  \\
\fs  \end{pmatrix}  = \sum \dfrac{\varepsilon_1^{\deg a_1} \dots \varepsilon_r^{\deg a_r }}{a_1^{s_1} \dots a_r^{s_r}}  \in K_{\infty}.
\end{equation*}
Here, the sum is over all tuples of monic polynomials with decreasing degree as MZV's, i.e., $(a_1,\ldots,a_r) \in A_+^r$ with $\deg (a_1)>\cdots>\deg (a_r)$. We define the depth and weight of $ \zeta_A \begin{pmatrix}
 \fve  \\
\fs  \end{pmatrix} $ in a similar manner, that is, $r$ is the depth and $w(\fs) = s_1 + \dots + s_r$ is the weight.

In \textit{loc. cit.}, Harada also showed some fundamental properties of AMZV's which hold for MZV's including non-vanishing, sum-shuffle relations, period interpretation and linear independence.

\subsubsection{Zagier-Hoffman's conjectures for AMZV's in positive characteristic} \ppar

The primary objective of the theory of MZV's and AMZV's in positive characteristic is to comprehend all of their $K$-linear relations as in the classical setting. The positive characteristic version of Zagier-Hoffman's conjectures for MZV's is suggested by Todd \cite{Tod18} and Thakur \cite{Tha17}. The algebraic part of these conjectures was solved by the fourth author \cite{ND21} using ingredients of Chen \cite{Che15}, Thakur \cite{Tha10, Tha17} and Todd \cite{Tod18}. The transcendental part was solved for the small weight cases by the fourth author in \cite{ND21} by applying Anderson's theory of motives \cite{And86,BP20,HJ20}, as well as Anderson-Brownawell-Papanikolas criterion \cite{ABP04} as a tool for transcendence.

In \cite{IKLNDP22}, by discovering a deep link among AMZV's and alternating Carlitz multiple polylogarithms, we solved both the MZV's and AMZV's version of the Zagier-Hoffman's conjectures in positive characteristic. We mention that the latter is much harder than the former (see \cite[\S 0.3]{IKLNDP22} for a detailed discussion). 

More precisely, we show that 
	\[ \dim_K \mathcal{AZ}_w = s(w) \]
where $\mathcal{AZ}_w$ is the vector space spanned by the AMZV's in the positive characteristic of weight $w$ over $K$, and $s(w)$ is a Fibonacci-like sequence defined by
\begin{align*}
s(w)=\begin{cases}
1 & \text{ if } w=0, \\
(q-1) q^{w-1} & \text{ if } 1 \leq w \leq q-1, \\
(q-1) (q^{w-1}-1) & \text{ if } w=q,\\
(q-1) \sum_{i=1}^{q-1} s(w-i)+s(w-q) & \text{ if } w>q.
\end{cases}
\end{align*}

Further, we give an explicit basis over $K$ of $\mathcal{AZ}_w$. Unlike the MZV's case, this basis consists of linear combinations of AMZV’s, not AMZV's. This crucial point has been observed in the classical setting by Charlton \cite{Cha21} and Deligne \cite{Del10}.

\subsubsection{Algebraic structures of MZV's in positive characteristic} \ppar

In our recent paper \cite{IKLNDP23} we addressed a question raised by a referee of a previous paper \cite{ND21} by the fourth author investigating the algebraic structure of MZV's in positive characteristic. It turns out\footnote{We are thankful to Thakur for his assistance with event orders.} that in 2016, Shuji Yamamoto \cite{Tha17} gave a construction of the word algebra for the shuffle product, and Shi \cite{Shi18} equipped it with a conjectural Hopf algebra structure, but left some axioms for the Hopf algebra open (see \cite[Conjecture 3.2.2 and  3.2.11]{Shi18}). Following this, Thakur wrote to Deligne, and he \cite{Del17} responded with a helpful letter in 2017. In it, Deligne stated that he had never seen such a construction before and advised on how to proceed with proving and investigating it.

\begin{conjecture}\label{conj: MZV Hopf}
    The set of MZV's in positive characteristic can be equipped with a shuffle Hopf algebra structure.
\end{conjecture}

In \cite{IKLNDP23} we proved Conjecture~\ref{conj: MZV Hopf}:
\begin{theorem}[{\cite[Theorem B]{IKLNDP23}}]
Conjecture~\ref{conj: MZV Hopf} is true.
\end{theorem}

Further, based on the previous link among AMZV's and alternating Carlitz multiple polylogarithms introduced in \cite{IKLNDP22}, we deduced a Hopf algebra of MZV's for the stuffle product (see \cite[\S 9]{IKLNDP23}).

\subsection{Main results and key ingredients}

\subsubsection{Main results} ${}$\par

In this paper we present a systematic study of algebraic structures of AMZV's in positive characteristic, and hence extend our previous work for MZV's. It could be considered as an analogue of Hoffman's and Racinet's results mentioned above in the classical setting.

Let us stress that the extension from the setting of MZV's to that of AMZV's is highly nontrivial. This fact is somewhat expected as in \cite{IKLNDP22} we have encountered a similar problem while extending the proof of Zagier-Hoffman's conjectures for MZV's to that for AMZV's.

Roughly speaking, in the case of MZV's, if we consider the corresponding word algebra denoted by $\mathfrak C$ (see \S \ref{sec: shuffle algebra MZV's}), then there is only one letter of fixed weight $n$ which is denoted by $x_n$. Consequently, it follows that there is also a unique and canonical way to define the coproducts $\Delta(x_n)$ so that the compatibility for $x_1$ and $x_n$ is satisfied for all $n \in \bN$ (see \S \ref{sec: coproduct for MZV's}). Unfortunately, when we work with AMZV's, the situation is much more complicated. Although it is not difficult to define the word algebra $\mathfrak D$ of AMZV's and to show that it is associative, the definition of a coproduct is far from evident. The main issue is that in $\mathfrak D$ there are $q-1$ letters of fixed weight, and there are no canonical ways to define coproducts for letters by using the compatibility for letters of weight $1$ and other letters. To move forward we would need some new ideas.

\subsubsection{Key ingredients} ${}$\par

The solution that we give here differs from that for the MZV's case at some crucial points and requires two new ingredients.

The first new ingredient is giving direct and explicit formulas for the coproduct $\Delta$ (see \S \ref{sec: coproduct}). We do not define it by the recursive induction on weight as it is the case for MZV's.

The second new ingredient is introducing so-called horizontal maps which allow us to deal with letters of the same weight (see \S \ref{sec: horizontal maps}). In particular, if some property is verified for a word beginning with a specific letter, says $x$, then by using horizontal maps we are able to deduce this property for all words beginning with a letter of the same weight as $x$.

With these tools we prove that the coproduct verifies all the desired properties so that the word algebra $\mathfrak D$ of AMZV's is indeed a connected graded Hopf algebra.

\subsubsection{Generalization for colored MZV's in positive characteristic} ${}$\par

We now give some comments about our results.

First, as mentioned above we have discovered a connection among AMZV's and the counterpart of Carlitz multiple polylogarithms in our previous work \cite[\S 1]{IKLNDP22}. Using this connection we could define another Hopf algebra structure on $\mathfrak D$ called the stuffle Hopf algebra of AMZV's which extends \cite[\S 9]{IKLNDP23}. As the extension can be done without difficulty, we leave it to the reader as an exercise.

Second, our construction of the Hopf algebra structure for the word algebra $\mathfrak D$ of AMZV's can be extended in a straightforward manner to the word algebra  $\mathfrak D_N$ of colored MZV's for all $N \in \mathbb N$. More precisely, the letters of $\mathfrak D_N$ consist of $x_{n,\varepsilon}$ where $n \in \mathbb N$ and $\varepsilon$ is an $N$-th root of unity in $\overline{\mathbb F}_q$. However, we warn the reader that in general, colored MZV's do not belong to $K_\infty$ but only to a finite extension of $K_\infty$.

\subsubsection{Plan of the manuscript} ${}$\par

The structure of this paper is as follows.
\begin{itemize}
\item In \S \ref{sec: review of Hopf algebra} we review some notions on the Hopf algebra. Our presentation is identical to that given in \cite{IKLNDP23}.

\item In \S \ref{sec: algebraMZVs}, we recall the word algebra $\mathfrak C$ of MZV's (called the composition space in \cite{IKLNDP23}) and the shuffle Hopf algebra structure of $\mathfrak C$. Then we review main results of {\it loc. cit} and state all necessary results that we will need in this paper so that the reader could read it easily without referring to \cite{IKLNDP23}.

\item In \S \ref{sec: review of AMZV's} we review the definitions and basic properties of the AMZV’s in positive characteristic introduced by Harada \cite{Har21}.

\item In \S \ref{sec: shuffle algebra} we introduce the word algebra or composition space $\mathfrak{D}$ of AMZV's and equip it with the shuffle algebra structure. We then define horizontal maps in \S \ref{sec: horizontal maps}. Finally, using these maps we prove that $\mathfrak D$ is an algebra (i.e., the associativity) and that there exists a shuffle map for AMZV's (see Theorem~\ref{thm: commutative algebra} and Theorem~\ref{thm: AZ shuffle map}).

\item In \S \ref{sec: shuffle Hopf algebra} we define the coproduct on $\mathfrak{D}$. This definition differs from that on $\mathfrak C$ as it is direct, explicit and does not use any induction on weight. Using horizontal maps we prove the compatibility and coassociativity results (see Theorem~\ref{thm: main statement on compatibility} and Theorem~\ref{thm: coassociativity}, respectively), establishing the main theorem of this paper (see {Theorem~\ref{thm: AMZV Hopf algebra}} for a precise statement):
\end{itemize}

\begin{theorem}
The word algebra $\mathfrak{D}$ of AMZV's equipped with the shuffle product and the above coproduct is a connected graded Hopf algebra of finite type over $\mathbb{F}_q$.
\end{theorem}

\subsection*{Acknowledgments}  

We would like to thank Dinesh Thakur for his interest and guidance with literature and for helpful discussions. Bo-Hae Im was supported by Basic Science Research Program through the National Research Foundation of Korea (NRF) funded by grant funded by the Korea government (MSIT) (NRF-2023R1A2C1002385). KN. Le and T. Ngo Dac were partially supported by the Excellence Research Chair ``$L$-functions in positive characteristic and applications'' financed by the Normandy Region.  T. Ngo Dac and LH. Pham were partially supported by the Vietnam Academy of Science and Technology (VAST) under grant no.~CTTH00.02/23-24 ``Arithmetic and Geometry of schemes over function fields and applications''.


\section{Review on Hopf algebra} \label{sec: review of Hopf algebra}

We give a brief review on basic notions of Hopf algebras, referring to our previous paper \cite[\S 2]{IKLNDP22} on the Hopf algebra structure of MZV's in positive characteristic.
\subsection{Hopf algebra} \ppar
Let $k$ be a field. From now on, tensor products are taken over $k$. Let $H$ be a $k$-vector space, and $\iota \colon H \otimes H \to H \otimes H$ be the transposition map, defined as $\iota(x\otimes y) = y \otimes x$ and to be linear.
\begin{definition}
    An algebra over $k$ is a triple $(H, m,u)$ consisting of a $k$-vector space $H$ together with $k$-linear maps  $m \colon H \otimes H \rightarrow H$ called the multiplication and $u \colon k \rightarrow H$ called the unit such that the following diagrams are commutative:
   \begin{enumerate}
       \item associativity
    \begin{equation*}
    \begin{tikzcd}
H\otimes H \otimes H \arrow[r, "m \otimes \text{id}"] \arrow[d, " \text{id} \otimes m "'] & H \otimes H \arrow[d, "m"] \\
H \otimes H \arrow[r, "m"]  & H
 \end{tikzcd}
 \end{equation*}
 \item unitary
\begin{equation*}
    \begin{tikzcd}
H \otimes k \arrow[rd] \arrow[r, "\text{id} \otimes u"] & H \otimes H \arrow[d, "m"] & k \otimes H \arrow[ld] \arrow[l, "u \otimes \text{id}"'] \\
& H & \end{tikzcd}
    \end{equation*}
    where the diagonal arrows are canonical isomorphisms.
   \end{enumerate}
   The algebra is said to be commutative if the following diagram is commutative:
    \begin{center}
        \begin{tikzcd}
H \otimes H \arrow[rd, "m"'] \arrow[rr, "\iota"] &   & H \otimes H \arrow[ld, "m"] \\                                                  & H &
\end{tikzcd}
    \end{center}
\end{definition}
Coalgebra over $k$ is dual to the algebra over $k$, obtained by reversing the arrows.
\begin{definition}
    A coalgebra over $k$ is a triple $(H, \Delta,\epsilon)$ consisting of a $k$-vector space $H$ together with $k$-linear maps  $\Delta \colon H  \rightarrow H \otimes H$ called the coproduct and $\epsilon \colon H \rightarrow k$ called the counit such that the following diagrams are commutative:
\begin{enumerate}[$(1)$]
    \item  coassociativity
    \begin{center}
       \begin{tikzcd}
H \otimes H\otimes H                              & \otimes H \arrow[l, "\Delta\otimes \text{id} "'] \\
H \otimes H \arrow[u, "\text{id} \otimes \Delta"] & {H} \arrow[u, "\Delta"'] \arrow[l, "\Delta"']
\end{tikzcd}
    \end{center}
 \item counitary
    \begin{center}
    \begin{tikzcd}
H \otimes k & H \otimes H \arrow[l, "\text{id} \otimes \epsilon "'] \arrow[r, "\epsilon \otimes \text{id}"] & k \otimes H \\
            & H \arrow[u, "\Delta"'] \arrow[ru] \arrow[lu]                                                  &
\end{tikzcd}
    \end{center}
    where the diagonal arrows are canonical isomorphisms.
    \end{enumerate}
    The coalgebra is said to be cocommutative if the following diagram is commutative:
    \begin{center}
        \begin{tikzcd}
H \otimes H \arrow[rr, "\iota"] &                                              & H \otimes H \\
                            & H \arrow[ru, "\Delta"'] \arrow[lu, "\Delta"] &
\end{tikzcd}
    \end{center}
\end{definition}

\begin{definition}
    A bialgebra over $k$ is a 5-tuple $(H,m, u, \Delta,\epsilon)$ consisting of an algebra $(H,m, u)$ over $k$  and a coalgebra $(H,\Delta,\epsilon)$ over $k$   which are compatible, i.e., the following diagrams are commutative:
\begin{enumerate}[$(1)$]
    \item product and coproduct
    \begin{center}
      \begin{tikzcd}
H \otimes H \arrow[d, "\Delta \otimes \Delta"'] \arrow[r, "m"]                  & H \arrow[r, "\Delta"] & H \otimes H                                                          \\
H \otimes H \otimes H \otimes H  \arrow[rr, "\text{id} \otimes \iota \otimes \text{id}"] &                       & H \otimes H \otimes H \otimes H  \arrow[u, "m \otimes m"']
\end{tikzcd}
    \end{center}
 \item unit and coproduct
    \begin{center}
    \begin{tikzcd}
H \arrow[r, "\Delta"]         & H \otimes H                                 \\
k \arrow[u, "u"] \arrow[r, "\text{canonical}"] & k \otimes k \arrow[u, "u \otimes u"']
\end{tikzcd}
    \end{center}
\item counit and product
    \begin{center}
    \begin{tikzcd}
H \arrow[d, "\epsilon"'] & H \otimes H \arrow[d, "\epsilon \otimes \epsilon"] \arrow[l, "m"'] \\
k                        & k \otimes k \arrow[l, "\text{canonical}"]
\end{tikzcd}
    \end{center}
  \item unit and counit
\begin{center}
    \begin{tikzcd}
k \arrow[rr, "\text{id}"] \arrow[rd, "u"'] &                           & k \\
                                              & H \arrow[ru, "\epsilon"'] &
\end{tikzcd}
\end{center}
\end{enumerate}
\end{definition}

\begin{definition}
    A Hopf algebra over $k$ is a bialgebra $(H,m, u, \Delta,\epsilon)$ over $k$ together with a $k$-linear map $S\colon H \rightarrow H$ called antipode such that the following diagram is commutative:
    \begin{equation*}
        \begin{tikzcd}
& H \otimes H \arrow[rr, "S \otimes \text{id}"] &                      & H \otimes H \arrow[rd, "m"]   &   \\
H \arrow[rr, "\epsilon"] \arrow[ru, "\Delta"] \arrow[rd, "\Delta"'] &                                               & k \arrow[rr, "u"] &                                    & H \\
& H \otimes H \arrow[rr, "\text{id}\otimes S "] &                      & H \otimes H \arrow[ru, "m"'] &
\end{tikzcd}
    \end{equation*}
\end{definition}

We note that there exist a bialgebra which does not admit an antipode (see \cite[Exercise 3.83]{BGF} for an example).

\subsection{Graded Hopf algebras} ${}$\par
We recall the notions of connected graded bialgebras.

\begin{definition} \label{defn: graded Hopf algebra}
\text{ }
\begin{enumerate}[$(1)$]
    \item   A bialgebra $(H,m, u, \Delta,\epsilon)$ over $k$ is said to be graded if one can write $H$ as a direct sum of $k$-subspaces
    \begin{equation*}
        H = \bigoplus \limits_{n = 0}^{\infty}H_n,
    \end{equation*}
    such that for all integers $r,s\geq 0$, we have
    \begin{equation*}
        m(H_r \otimes H_s) \subseteq H_{r + s} \quad \text{and} \quad \Delta(H_r) \subseteq \bigoplus \limits_{i + j = r} H_i \otimes H_j.
    \end{equation*}
    A graded bialgebra is said to be connected if $H_0 = k$.
\item A graded Hopf algebra is a Hopf algebra $H$ whose underlying bialgebra is graded and the antipode $S$ satisfies $S(H_n) \subseteq H_n$.

\item A graded Hopf algebra is said to be connected if $H_0 = k$.

\item A graded Hopf algebra is said to be {\it of finite type} if $H_n$ is a $k$-vector space of finite dimension.
\end{enumerate}
\end{definition}

 The following proposition shows that a connected graded bialgebra automatically admits an antipode, thus it is always a Hopf algebra.
\begin{proposition} \label{prop: graded Hopf algebras}
    Let $(H,m, u, \Delta,\epsilon)$ be a connected graded bialgebra over $k$.
\begin{enumerate}[\normalfont (1)]
\item For each element $x \in H_n$ with $n \geq 1$, we have
    \begin{equation*}
         \Delta(x) = 1 \otimes x + x \otimes 1 + \sum x_{(1)} \otimes x_{(2)},
    \end{equation*}
    where $\sum   x_{(1)} \otimes x_{(2)} \in \bigoplus \limits_{\substack{i,j > 0 \\ i + j = n}} H_i \otimes H_j$. Moreover, the counit $\epsilon$ vanishes on~$H_n$ for all $n \geq 1$.
\item  We continue the notation as in $(1)$ and define recursively a $k$-linear map $S \colon H \rightarrow H$ given by
    \begin{align*}
S(x) = \begin{cases}  x & \quad \text{if } x \in H_0, \\
- x - \sum  m(S(x_{(1)}) \otimes x_{(2)}) & \quad \text{if } x \in H_n \text{ with } n \geq 1.
\end{cases}
\end{align*}
Then $H$ is a graded Hopf algebra whose antipode is $S$.
\end{enumerate}
\end{proposition}

\begin{proof}
See \cite[Lemma 2.1]{Ehr96}.
\end{proof}


\section{Review of shuffle Hopf algebra of MZV's} \label{sec: algebraMZVs}

We briefly review some results of our previous work  \cite{IKLNDP23} on algebra structures of MZV's in positive characteristic.

\subsection{Alphabet} \label{sec: alphabet} ${}$\par

Let $I$ be a countable set and $A = \{x_n\}_{n \in I}$ be a set of variables indexed by $I$, equipped with  weights $w(x_n) \in \mathbb N$. The set $A$ will be called an alphabet and its elements will be called letters. A word over the alphabet $A$ is a finite string of letters. In particular, the empty word will be denoted by $1$. The depth $\depth(\fa)$ of a word $\fa$ is the number of letters in the string of $\fa$, so that $\depth(1) = 0$. The weight of a word is the sum of the weights of its letter and we put $w(1)=0$. Let $\langle A \rangle$ denote the set of all words over $A$. We endow $\langle A \rangle$ with the concatenation product defined by the following formula:
\begin{equation*}
    x_{i_1} \dotsc x_{i_n} \cdot x_{j_1} \dotsc x_{j_m} = x_{i_1} \dotsc x_{i_n} x_{j_1} \dotsc x_{j_m}.
\end{equation*}
Let $\Fq \langle A \rangle$ be the free $\Fq$-vector space with basis $\langle A \rangle$. The concatenation product extends to $\Fq \langle A \rangle$ by linearity. For a letter $x_a \in A$ and an element $\fa \in \Fq \langle A \rangle$, we write simply $x_a\fa$ instead of $x_a \cdot \fa$. For each nonempty word $\fa \in \langle A \rangle$, we can write $\fa = x_a \fa_-$ where $x_a$ is the first letter of $\fa$ and $\fa_-$ is the word obtained from $\fa$ by removing $x_a$.

\subsection{Alphabet and shuffle algebra associated to MZV's} \label{sec: shuffle algebra MZV's} ${}$\par

We recall the alphabet associated to MZV's as given in \cite[\S 4.2]{IKLNDP23}. The set $I$ will be $\mathbb N$ and we denote by $\Sigma = \{x_n\}_{n \in \mathbb{N}}$ with weight $w(x_n)=n$ the alphabet attached to MZV's. We put $\mathfrak{C}=\Fq \langle \Sigma \rangle$, i.e., the free $\Fq$-vector space with basis $\langle \Sigma \rangle$.

For positive integers $r,s$ and for all positive integers $i,j$ such that $i+j=r+s$, we put
\begin{align} \label{eq: delta}
\Delta_{r,s}^{i}=\begin{cases}  (-1)^{r-1}\binom{i - 1}{r - 1} + (-1)^{s-1}\binom{i - 1}{s - 1} & \quad \text{if } (q - 1) \mid i, \\
0 & \quad \text{otherwise}.
\end{cases}
\end{align}
Then we define recursively two products on $\mathfrak{C}$ as $\Fq$-bilinear maps
\begin{align*}
   \diamond \colon \mathfrak{C} \times \mathfrak{C} \longrightarrow \mathfrak{C} \quad \text{and} \quad
   \shuffle \colon \mathfrak{C} \times \mathfrak{C} \longrightarrow \mathfrak{C}
\end{align*}
by setting $1 \diamond \mathfrak{a} = \mathfrak{a} \diamond 1 = \mathfrak{a}, 1 \shuffle  \mathfrak{a} = \mathfrak{a} \shuffle  1 = \mathfrak{a}$ and
\begin{align*}
    \fa \diamond \fb =& \ x_{a + b}(\fa_- \shuffle  \fb_-) + \sum\limits_{i+j = a + b} \Delta^j_{a,b} x_i(x_j \shuffle  (\fa_- \shuffle  \fb_-)),\\
    \fa \shuffle  \fb =& \ x_{a}(\fa_- \shuffle  \fb) + x_{b}(\fa \shuffle  \fb_-) + \fa \diamond \fb
\end{align*}
for any nonempty words $\fa,\fb \in \langle \Gamma \rangle$. We call $\diamond$ the diamond product and $\shuffle $ the shuffle product. The unit $u:\Fq \to \frak C$ is given by sending $1$ to the empty word. One of the main results of \cite{IKLNDP23} reads as follows (see \cite[Theorem A]{IKLNDP23}):

\begin{theorem} \label{thm: shuffle algebra}
The spaces $(\mathfrak{C}, \diamond)$ and $(\mathfrak{C}, \shuffle )$ are commutative $\Fq$-algebras.
\end{theorem}

\subsection{Shuffle Hopf algebra associated to MZV's} \label{sec: coproduct for MZV's} ${}$\par

We briefly recall the coproduct as in \cite[\S 7.1]{IKLNDP23}
	\[ \Delta: \frak C \to \frak C \otimes \frak C. \]
We first define it on  $\langle \Sigma \rangle$ by induction on weight and extend by $\Fq$-linearity to $\frak C$. We put
\begin{align*}
\Delta(1):=1 \otimes 1, \quad\quad\Delta(x_1):=1 \otimes x_1 + x_1 \otimes 1.
\end{align*}
Let $w \in \bN$ and we suppose that we have defined $\Delta(\fu)$ for all words $\fu$ of weight $w(\fu)<w$. We now give a formula for $\Delta(\fu)$ for all words $\fu$ with $w(\fu)=w$. For such a word $\fu$ with $\depth(\fu)>1$, we put $\fu=x_u \fv$ with $w(\fv)<w$. Since $x_u$ and $\fv$ are both of weight less than $w$, we have already defined
\begin{align*}
\Delta(x_u)=1 \otimes x_u + \sum a_u \otimes b_u, \quad\quad
\Delta(\fv)= \sum a_\fv \otimes b_\fv.
\end{align*}
Then we set
\begin{align*}
\Delta(\fu):=1 \otimes \fu + \sum (a_u \triangleright a_\fv) \otimes (b_u \shuffle  b_\fv).
\end{align*}
Our last task is to define $\Delta(x_w)$. We know that
	\[ x_1\shuffle x_{w-1}=x_w+x_1x_{w-1}+x_{w-1}x_1+\sum_{0<j<w} \Delta^j_{1,w-1}  x_{w-j} x_j \]
where all the words $x_{w-j} x_j$ have weight $w$ and depth $2$ and all $\Delta^j_{1,w-1}$ defined as in \eqref{eq: delta} belong to~$\Fq$. Therefore, we set
\begin{equation*}
\Delta(x_w):=\Delta(x_1) \shuffle  \Delta(x_{w-1})-\Delta(x_1x_{w-1})-\Delta(x_{w-1}x_1)-\sum_{0<j<w} \Delta^j_{1,w-1} \Delta(x_{w-j}x_j).
\end{equation*}
By induction on weight, we can check that $\Delta$ preserves the grading.

Next, the counit $\epsilon:\frak C \to \Fq$ is defined by putting $\epsilon(1)=1$ and $\epsilon(\fu)=0$ otherwise. We proved (see \cite[Theorem B]{IKLNDP23}):

\begin{theorem} \label{thm: MZV hopf algebra}
The connected graded bialgebra $(\frak C,\shuffle,u,\Delta,\epsilon)$ is a connected graded Hopf algebra over $\Fq$.
\end{theorem}

\subsection{Some properties of the coproduct} ${}$\par

In this section we collect some properties of the coproduct that will be useful in the sequel.

\begin{lemma} \label{lem: au non empty}
For all words $\fu \in \frak C$, we have
	\[ \Delta(\fu)=1 \otimes \fu+\sum a_\fu \otimes b_\fu, \text{ where } a_\fu \neq 1. \]
\end{lemma}

\begin{proof}
See \cite[Lemma 7.1]{IKLNDP23}.
\end{proof}

\begin{lemma} \label{lem: delta diamond}
Let $\mathfrak{u}$ and $\mathfrak{v}$ be two nonempty words in $\langle \Sigma \rangle$. Suppose that
\begin{align*}
    \Delta(\mathfrak{u}) = 1 \otimes \mathfrak{u} + \sum \mathfrak{u}_{(1)} \otimes \mathfrak{u}_{(2)}, \text{ and } \quad \Delta(\mathfrak{v}) = 1 \otimes \mathfrak{v} + \sum \mathfrak{v}_{(1)} \otimes \mathfrak{v}_{(2)}.
\end{align*}
Then
\begin{equation*}
    \Delta(\mathfrak{u} \diamond \mathfrak{v}) = 1 \otimes (\mathfrak{u} \diamond \mathfrak{v}) + \sum (\mathfrak{u}_{(1)} \diamond \mathfrak{v}_{(1)}) \otimes (\mathfrak{u}_{(2)} \shuffle \mathfrak{v}_{(2)}).
\end{equation*}
\end{lemma}

\begin{proof}
See \cite[Lemma 8.5]{IKLNDP23}.
\end{proof}

Finally, there is a useful formula for $\Delta(x_n)$ for $x_n \in \Sigma$ that we briefly explain below.

\begin{lemma} \label{lem: delta j<n}
Let $n \in \N$. Then we have
\begin{enumerate}[\normalfont (1)]
    \item for all $j < n$,
\begin{equation*}
    \Delta^j_{1,n} = \begin{cases}  1 & \quad \text{if } (q - 1) \mid j \\
0 & \quad \text{otherwise}.
\end{cases}
\end{equation*}

\item $\Delta^n_{1,n} = 0$.
\end{enumerate}
\end{lemma}

\begin{proof}
The result is straightforward from the definition of $\Delta^j_{1,n}$. See \cite[Lemma~8.1]{IKLNDP23}.
\end{proof}

Next, we introduce the bracket operator by the following formula: for any $\fa=x_{i_1} \dots x_{i_m} \in \langle \Sigma \rangle$,
\begin{equation} \label{eq: bracket operator}
    [\fa] := (-1)^{m}\Delta^{i_1}_{1,w(\fa) + 1}\dotsb\Delta^{i_m}_{1,w(\fa) + 1} x_{i_1} \shuffle \dotsb \shuffle x_{i_m}.
\end{equation}
As a matter of convention, we also agree that $[1] = 1$. Then we showed (see \cite[Theorem C]{IKLNDP23}):

\begin{proposition} \label{prop: formula delta xn}
We keep the above notation. Then for all $n \in \N$, we have
\begin{equation*}
    \Delta(x_n) = 1 \otimes x_n + \sum \limits_{\substack{r \in \N, \fa \in \langle \Sigma \rangle\\ r + w(\fa) = n}} \binom{r + \depth(\fa) - 2}{\depth(\fa)} x_r \otimes [\fa].
\end{equation*}
\end{proposition}


\section{Alternating multiple zeta values} \label{sec: review of AMZV's}
We recall some notions on MZV's and AMZV's in positive characteristic, referring to \cite{Har21, IKLNDP22, IKLNDP23}.

Letting $\fs = (s_1, \dotsc, s_n) \in \mathbb{N}^{n}$ and $\fve = (\varepsilon_1, \dotsc, \varepsilon_n) \in (\mathbb{F}_q^{\times})^{n}$, we set $\fs_- := (s_2, \dotsc, s_n)$ and $\fve_- := (\varepsilon_2, \dotsc, \varepsilon_n) $. A positive array $\begin{pmatrix}
 \fve  \\
\fs  \end{pmatrix} $ is an array of the form $$\begin{pmatrix}
 \fve  \\
\fs  \end{pmatrix}  = \begin{pmatrix}
 \varepsilon_1 & \dotsb & \varepsilon_n \\
s_1 & \dotsb & s_n \end{pmatrix}.$$

We recall the power sums studied by Thakur \cite{Tha10}. For $d \in \mathbb{Z}$ and for $\fs=(s_1,\dots,s_n) \in \N^n$ we introduce
\begin{equation*}
    S_d(\fs) =  \sum\limits_{\substack{a_1, \dots, a_n \in A_{+} \\ d = \deg a_1> \dots > \deg a_n\geq 0}} \dfrac{1}{a_1^{s_1} \dots a_n^{s_n}} \in K
\end{equation*}
and
\begin{equation*}
    S_{<d}(\mathfrak s) = \sum\limits_{\substack{a_1, \dots, a_n \in A_{+} \\ d > \deg a_1> \dots > \deg a_n\geq 0}} \dfrac{1}{a_1^{s_1} \dots a_n^{s_n}} \in K.
\end{equation*}
We also recall the following analogues of the power sums after Harada \cite{Har21}. For a positive array $\begin{pmatrix}
 \fve  \\
\fs  \end{pmatrix}  =  \begin{pmatrix}
 \varepsilon_1 & \dots & \varepsilon_n \\
s_1 & \dots & s_n \end{pmatrix} $, we introduce
\begin{equation*}
S_d \begin{pmatrix}
\fve \\ \fs
\end{pmatrix}  = \sum\limits_{\substack{a_1, \dots, a_n \in A_{+} \\ d = \deg a_1> \dots > \deg a_n\geq 0}} \dfrac{\varepsilon_1^{\deg a_1} \dots \varepsilon_n^{\deg a_n }}{a_1^{s_1} \dots a_n^{s_n}} \in K
\end{equation*}
and
\begin{equation*}
    S_{<d} \begin{pmatrix}
 \fve  \\
\fs  \end{pmatrix}  = \sum\limits_{\substack{a_1, \dots, a_n \in A_{+} \\ d > \deg a_1> \dots > \deg a_n\geq 0}} \dfrac{\varepsilon_1^{\deg a_1} \dots \varepsilon_n^{\deg a_n }}{a_1^{s_1} \dots a_n^{s_n}} \in K.
\end{equation*}
One verifies easily the following formulas:
\begin{align*}
    S_{d} \begin{pmatrix}
 \varepsilon \\
s  \end{pmatrix}  &= \varepsilon^d S_d(s), \\
S_d \begin{pmatrix}
 1& \dots & 1 \\
s_1 & \dots & s_n \end{pmatrix}  &= S_{d}(s_1, \dots, s_n),\\
S_{<d} \begin{pmatrix}
 1& \dots & 1 \\
s_1 & \dots & s_n \end{pmatrix}  &= S_{<d}(s_1, \dots, s_n), \\
S_{d} \begin{pmatrix}
 \fve  \\
\fs  \end{pmatrix}  &= S_{d} \begin{pmatrix}
 \varepsilon_1  \\
s_1  \end{pmatrix} S_{<d} \begin{pmatrix}
 \fve_{-}  \\
\fs_{-}  \end{pmatrix}.
\end{align*}

In \cite{Har21} Harada introduced the alternating multiple zeta value (AMZV) in positive characteristic by setting
\begin{equation*}
    \zeta_A \begin{pmatrix}
 \fve  \\
\fs  \end{pmatrix}  := \sum \limits_{d \geq 0} S_d \begin{pmatrix}
 \fve  \\
\fs  \end{pmatrix}  = \sum\limits_{\substack{a_1, \dots, a_n \in A_{+} \\ \deg a_1> \dots > \deg a_n\geq 0}} \dfrac{\varepsilon_1^{\deg a_1} \dots \varepsilon_n^{\deg a_n }}{a_1^{s_1} \dots a_n^{s_n}}  \in K_{\infty}.
\end{equation*}
Then he proved several fundamental properties of these values. In particular, it is shown that the product of two AMZV's can be expressed as a linear combination with coefficients in $\Fq$ of AMZV's. To do so, we first recall Chen's formula in \cite{Che15}. For integers $a, b$ with $b \geq 0
$, we recall the binomial number defined by
\begin{align*}
\binom{a}{b} :=& \dfrac{a(a-1) \dotsc (a-b+1)}{b!}.
\end{align*}
It should be remarked that $\binom{a}{b} = 0$ if $b > a \geq 0$. Then refining the work of Thakur~\cite{Tha10}, Chen \cite[Theorem 3.1]{Che15} showed that for positive integers $r,s$ and for all $d \in \mathbb{N}$,
\begin{equation*}
    S_d(r) S_d(s) = S_d(r+s) + \sum \limits_{i + j = r + s} \Delta_{r,s}^{j} S_d(i,j),
\end{equation*}
where the coefficients $\Delta_{r,s}^{j}$ are defined as in \eqref{eq: delta}.

As a  direct consequence of Chen's formula, we obtain

\begin{lemma} \label{lemma: VMPLdepth1}
Let $\alpha, \beta$ be two elements in $\mathbb{F}_q^*$, and let $r,s$ be two positive integers. Then for all $d \in \mathbb{N}$, we have
\begin{align*}
    S_d  \begin{pmatrix}
 \alpha \\
r  \end{pmatrix}  S_d \begin{pmatrix}
 \beta \\
s  \end{pmatrix} = S_d \begin{pmatrix}
 \alpha\beta \\
r + s  \end{pmatrix}
 + \sum \limits_{i + j = r + s} \Delta_{r,s}^{j} S_d \begin{pmatrix}
 \alpha\beta &   1\\
i & j  \end{pmatrix},
\end{align*}
where the indices $i,j$ are positive integers.
\end{lemma}

\begin{proof}
See for example \cite[Lemma 2.3]{Har21}.
\end{proof}


\section{Shuffle algebra of AMZV's} \label{sec: shuffle algebra}

In the following sections \S \ref{sec: shuffle algebra} and \S \ref{sec: shuffle Hopf algebra} we construct the Hopf shuffle algebra of AMZV's and note that the extension from the setting of MZV's to that of AMZV's is not straightforward as we have explained in the introduction. The key ingredient of our construction is the notion of so-called horizontal maps that we introduce in this section. These maps will allow us to deal with letters of the same weight (see \S \ref{sec: horizontal maps}). In the following section, we will give direct and explicit formulas for the coproduct and then prove that such formulas verify all desired properties.

\subsection{Alphabet and shuffle algebra associated to AMZV's} ${}$\par

We put $I = \{(n, \varepsilon) : n \in \mathbb{N}, \varepsilon \in \F_q^*\}$. Let $\Gamma = \{x_{n,\varepsilon}\}_{n \in \mathbb{N}, \varepsilon \in \mathbb{F}_q^*}$ be the alphabet associated to AMZV's indexed by $I$ (see \S \ref{sec: alphabet} for notation). We define the weights by $w(x_{n,\varepsilon})=n$. We define $\mathfrak{D}:= \Fq \langle \Gamma \rangle$ to be the free $\mathbb{F}_q$-vector space with basis~$\langle \Gamma \rangle$.

Inspired by \S \ref{sec: algebraMZVs} we define recursively two products on $\mathfrak{D}$ as $\mathbb{F}_q$-bilinear maps
\begin{align*}
   \diamond \colon \mathfrak{D} \times \mathfrak{D} \longrightarrow \mathfrak{D} \quad \text{and} \quad
   \shuffle \colon \mathfrak{D} \times \mathfrak{D} \longrightarrow \mathfrak{D}
\end{align*}
by setting $1 \diamond \mathfrak{a} = \mathfrak{a} \diamond 1 = \mathfrak{a},$ $ 1 \shuffle \mathfrak{a} = \mathfrak{a} \shuffle 1 = \mathfrak{a}$ and
\begin{align*}
    \fa \diamond \fb =& \ x_{a + b, \alpha \beta}(\fa_- \shuffle \fb_-) + \sum\limits_{i+j = a + b} \Delta^j_{a,b} x_{i, \alpha \beta}(x_{j, 1} \shuffle (\fa_- \shuffle \fb_-)),\\
    \fa \shuffle \fb =& \ x_{a,\alpha}(\fa_- \shuffle \fb) + x_{b,\beta}(\fa \shuffle \fb_-) + \fa \diamond \fb
\end{align*}
for any nonempty words $\fa ,\fb \in \langle \Gamma \rangle$ with $\fa = x_{a,\alpha} \fa_-,\fb = x_{b,\beta} \fb_- $. We call $\diamond$ the diamond product and $\shuffle $ the shuffle product.

\begin{proposition} \label{prop: commutative}
 The diamond product and the shuffle product on $\mathfrak{D}$ are commutative.
\end{proposition}

\begin{proof}
Let $\fa, \fb \in \langle \Gamma \rangle$ be two arbitrary words. It suffices to show that
\begin{equation*} \label{eq: commutative}
    \fa \diamond \fb = \fb \diamond \fa \quad \text{and} \quad \fa \shuffle \fb = \fb \shuffle \fa .
\end{equation*}
We proceed the proof by induction on $\text{depth}(\fa) + \text{depth}(\fb)$. If one of $\fa$ or $\fb$ is the empty word, then \eqref{eq: commutative} holds trivially. We assume that \eqref{eq: commutative} holds when $\text{depth}(\fa) + \text{depth}(\fb) < n$ with $n \in \mathbb{N}$ and $n \geq 2$. We need to show that \eqref{eq: commutative} holds when $\text{depth}(\fa) + \text{depth}(\fb) = n$.

Indeed, assume that $\fa = x_{a,\alpha} \fa_-,\fb = x_{b,\beta} \fb_- $, we have
\begin{align*}
    \fa \diamond \fb =& \ x_{a + b, \alpha \beta}(\fa_- \shuffle \fb_-) + \sum\limits_{i+j = a + b} \Delta^j_{a,b} x_{i, \alpha \beta}(x_{j, 1} \shuffle (\fa_- \shuffle \fb_-)),\\
    \fb \diamond \fa =& \ x_{b + a,  \beta\alpha}(\fb_- \shuffle \fa_-) + \sum\limits_{i+j = b + a} \Delta^j_{b,a} x_{i,  \beta\alpha}(x_{j, 1} \shuffle (\fb_- \shuffle \fa_-)).
\end{align*}
It follows from the induction hypothesis that $\fa_-\shuffle\fb_- = \fb_- \shuffle\fa_-$, hence $\fa \diamond \fb = \fb \diamond \fa$. On the other hand, we have
\begin{align*}
    \fa \shuffle \fb =& \ x_{a,\alpha}(\fa_- \shuffle \fb) + x_{b,\beta}(\fa \shuffle \fb_-) + \fa \diamond \fb,\\
    \fb \shuffle \fa =& \ x_{b,\beta}(\fb_- \shuffle \fa) + x_{a,\alpha}(\fb \shuffle \fa_-) + \fb \diamond \fa.
\end{align*}
It follows from the induction hypothesis and the above arguments that $\fa_- \shuffle \fb = \fb \shuffle \fa_-, \fa \shuffle \fb_- = \fb_- \shuffle \fa$ and $\fa \diamond \fb = \fb \diamond \fa$, hence $\fa \shuffle \fb = \fb \shuffle \fa$. This proves the proposition.
\end{proof}

\subsection{Horizontal maps} \label{sec: horizontal maps} ${}$\par

In this section, we introduce horizontal maps which are crucial in the sequel. For each element  $\alpha \in \mathbb{F}_q^*$, we consider the $\Fq$-linear map
\begin{equation*}
     \varphi_{\alpha} \colon \mathfrak{C} \longrightarrow \mathfrak{D},
\end{equation*}
which maps the empty word $1 \in \langle \Sigma \rangle$ to the empty word $1 \in \langle \Gamma \rangle$ and maps a nonempty word $x_{i_1}x_{i_2}\dotsc x_{i_n} $ to $ x_{i_1,\alpha}x_{i_2,1}\dotsc x_{i_n, 1}$. For all nonempty words $\fa \in \langle \Sigma \rangle$ with $\fa= x_a\fa_-$, we get
\begin{equation*}
    \varphi_{\alpha}(\fa) = \varphi_{\alpha}(x_a) \varphi_{1}(\fa_-) = x_{a, \alpha}\varphi_{1}(\fa_-).
\end{equation*}
As $w(\varphi_\alpha(\fa))=w(\fa)$ for all $\alpha \in \Fq^*$ and $\fa \in \frak C$, we call these maps horizontal maps. We note that they are all injective.

\begin{lemma} \label{lem: empty}
Let $\fa,\fb$ be two words in $\langle \Sigma \rangle$. Then we have
\begin{equation*} \label{eq: empty}
    \varphi_{1}(\fa \diamond \fb) = \varphi_{1}(\fa) \diamond \varphi_{1}(\fb) \quad \text{and} \quad \varphi_{1}(\fa \shuffle \fb) = \varphi_{1}(\fa) \shuffle \varphi_{1}(\fb).
\end{equation*}
\end{lemma}

\begin{proof}
The proof can be done by induction on $\text{depth}(\fa) + \text{depth}(\fb)$. We omit the details.
\end{proof}

As a direct consequence, $\varphi_1$ is a homomorphism of algebras from $(\frak C, \diamond)$ to $(\frak D, \diamond)$ (resp. from $(\frak C, \shuffle)$ to $(\frak D, \shuffle)$). 

\begin{notation}
Since $\varphi_1$ is injective, from now on for $\fa \in \frak C$ we write $\fa \in \frak D$ instead of $\varphi_1(\fa) \in \frak D$.
\end{notation}

Further, we extend the horizontal map as follows:
\[\varphi_\alpha \colon \frak D \to \frak D; \quad 1\mapsto 1, \quad \fu \mapsto x_{u, \alpha \varepsilon} \fu_-\]
for nonempty words $\fu \in \langle \Gamma \rangle$ with $\fu = x_{u, \varepsilon} \fu_-$.
This is a natural extension under the identification of $\fa \in \mathfrak{C}$ and $\varphi_1(\fa)\in \mathfrak{D}$. We note that if $\fu \in \langle \Gamma \rangle$ is a nonempty word written as $\fu = x_{u, \varepsilon} \fu_-$, then for all $\alpha \in \Fq^*$,
	\[ \varphi_\alpha(\fu) = \varphi_\alpha(x_{u, \varepsilon}) \fu_-. \]
Under this extension, the above lemma holds for $\fa, \fb \in \langle \Gamma \rangle$.

\begin{lemma} \label{lem: distributive}
Let $\alpha, \beta$ be elements of $\mathbb{F}_q^*$, and let $\fa,\fb$ be nonempty words in $\langle \Sigma \rangle$. We have
\begin{equation*}
   \varphi_{\alpha}(\fa) \diamond \varphi_{\beta}(\fb) = \varphi_{\alpha\beta}(\fa \diamond \fb).
\end{equation*}
Further, for nonempty words $\fa, \fb \in \langle \Gamma \rangle$, \begin{equation*}
   \varphi_{\alpha}(\fa) \diamond \varphi_{\beta}(\fb) = \varphi_{\alpha\beta}(\fa \diamond \fb).
\end{equation*}
\end{lemma}
\begin{proof}
Consider the horizontal maps defined on $\langle \Sigma \rangle$.
We write $\fa=x_a \fa_-, \fb=x_b \fb_-$ with $x_a, x_b \in \Sigma$ and $\fa_-,\fb_- \in \langle \Sigma \rangle$. From Lemma \ref{lem: empty}, we have
\begin{align*}
    \varphi_{\alpha\beta}(\fa \diamond \fb) =&  \ \varphi_{\alpha\beta}\Bigg(x_{a + b}(\fa_- \shuffle \fb_-) + \sum\limits_{i+j = a + b} \Delta^j_{a,b} x_i(x_j \shuffle (\fa_- \shuffle \fb_-))\Bigg)\\
    =& \ x_{a + b,\alpha\beta} (\fa_-\shuffle \fb_-) + \sum\limits_{i+j = a + b} \Delta^j_{a,b} x_{i,\alpha\beta}(x_j \shuffle (\fa_- \shuffle \fb_-))\\
    =& \ x_{a,\alpha}\fa_- \diamond x_{b,\beta}\fb_- \\
    =& \ \varphi_{\alpha}(\fa) \diamond \varphi_{\beta}(\fb).
\end{align*}

The proof for the horizontal maps on $\langle \Gamma \rangle$ is parallel. Let $\fa = x_{a, \varepsilon_1}\fa_-$, $\fb = x_{b, \varepsilon_2}\fb_-$ with $x_{a, \varepsilon_1}, x_{b, \varepsilon_2} \in \Gamma$ and $\fa_-,\fb_- \in \langle \Gamma \rangle$. From Lemma \ref{lem: empty} and its following remark, we have
\begin{align*}
    \varphi_{\alpha\beta}(\fa \diamond \fb) =&  \ \varphi_{\alpha\beta}\Bigg(x_{a + b, \varepsilon_1 \varepsilon_2}(\fa_- \shuffle \fb_-) + \sum\limits_{i+j = a + b} \Delta^j_{a,b} x_{i, \varepsilon_1 \varepsilon_2}(x_{j} \shuffle (\fa_- \shuffle \fb_-))\Bigg)\\
    =& \ x_{a + b,\alpha\beta \varepsilon_1 \varepsilon_2}(\fa_-\shuffle \fb_-) + \sum\limits_{i+j = a + b} \Delta^j_{a,b} x_{i,\alpha\beta \varepsilon_1 \varepsilon_2}(x_{j} \shuffle (\fa_- \shuffle \fb_-))\\
    =& \ x_{a,\alpha\varepsilon_1}\fa_- \diamond x_{b,\beta \varepsilon_2}\fb_-\\
    =& \ \varphi_{\alpha}(\fa) \diamond \varphi_{\beta}(\fb).
\end{align*}
This proves the lemma.
\end{proof}

We end this section by noting that horizontal maps are associative. More precisely, we prove
\begin{lemma} \label{lemma: associativity horizontal maps}
Let $\alpha, \beta \in \mathbb{F}_q^*$. For any word $\fa \in \langle \Gamma \rangle$, we have
\begin{equation*}
   \varphi_{\alpha \beta}(\fa) = \varphi_{\alpha} (\varphi_\beta(\fa)).
\end{equation*}
\end{lemma}

\begin{proof}
The proof is straightforward from the definition of horizontal maps.
\end{proof}

\subsection{Associativity} ${}$\par

In this subsection, we first show that the diamond product is associative over $\Gamma$ (see Proposition \ref{prop: associativiy of letters}), and then prove the associativity of the diamond product and the shuffle product on $\mathfrak{D}$ (see Proposition \ref{prop: associative}).

\begin{proposition} \label{prop: associativiy of letters}
Let $x_{a,\alpha},x_{b,\beta},x_{c,\gamma}$ be letters in $\Gamma$. Then we have
\begin{equation*}
    (x_{a,\alpha} \diamond x_{b,\beta}) \diamond x_{c,\gamma} = x_{a,\alpha} \diamond (x_{b,\beta} \diamond x_{c,\gamma}).
\end{equation*}
\end{proposition}

\begin{proof}
From Lemma \ref{lem: distributive}, we have
\begin{align*}
    (x_{a,\alpha} \diamond x_{b,\beta}) \diamond x_{c,\gamma} &= (\varphi_{\alpha}(x_{a}) \diamond \varphi_{\beta}(x_{b})) \diamond \varphi_{\gamma}(x_{c}) \\
    &=  \varphi_{\alpha\beta}(x_{a}\diamond x_{b}) \diamond \varphi_{\gamma}(x_{c})\\
    &= \varphi_{\alpha\beta\gamma}((x_a \diamond x_b)\diamond x_c),\\
    x_{a,\alpha} \diamond (x_{b,\beta} \diamond x_{c,\gamma}) &= \varphi_{\alpha}(x_{a}) \diamond (\varphi_{\beta}(x_{b}) \diamond \varphi_{\gamma}(x_{c})) \\
    &=   \varphi_{\alpha}(x_{a}) \diamond \varphi_{\beta\gamma}(x_{b}\diamond x_{c})\\
    &= \varphi_{\alpha\beta\gamma}(x_a \diamond (x_b\diamond x_c)).
\end{align*}
From the associativity of the diamond product on $\mathfrak{C}$ (see Theorem \ref{thm: shuffle algebra}), we have $(x_a \diamond x_b)\diamond x_c = x_a \diamond (x_b\diamond x_c)$. This proves the proposition.
\end{proof}

Next, we define recursively a product on $\mathfrak{D}$ as a $\Fq$-bilinear map
\begin{align*}
   \triangleright \colon \mathfrak{D} \times \mathfrak{D} \longrightarrow \mathfrak{D}
\end{align*}
by setting $1 \triangleright \mathfrak{a} = \mathfrak{a} \triangleright 1 = \mathfrak{a}$ and
\begin{align*}
    \fa \triangleright \fb =& \ x_{a, \alpha}(\fa_- \shuffle  \fb)
\end{align*}
for any nonempty words $\fa,\fb \in \langle \Gamma \rangle$ with $\fa = x_{a, \alpha}\fa_-$. We call $\triangleright$ the triangle product. The triangle product is neither commutative nor associative, as one verifies at once.

\begin{lemma} \label{lem: triangleright formulas}
For all words $\fa,\fb \in \langle \Gamma \rangle$ with $\fa = x_{a, \alpha}\fa_-, \fb = x_{b, \beta}\fb_-$, we have
\begin{enumerate}[\normalfont (1)]
    \item $\fa \diamond \fb = (x_{a, \alpha} \diamond x_{b, \beta}) \triangleright (\fa_- \shuffle  \fb_-)$.
    \item $\fa \shuffle  \fb = \fa \triangleright \fb + \fb \triangleright \fa + \fa \diamond \fb$.
\end{enumerate}
\end{lemma}

\begin{proof}
 We have
\begin{align*}
    \fa \diamond \fb
    =& \ x_{a + b, \alpha \beta}(\fa_- \shuffle \fb_-) + \sum\limits_{i+j = a + b} \Delta^j_{a,b} x_{i, \alpha \beta}(x_{j} \shuffle (\fa_- \shuffle \fb_-))\\
    =& \ x_{a + b, \alpha \beta} \triangleright (\fa_- \shuffle  \fb_-) + \sum \limits_{i+j = a + b} \Delta^j_{a,b} (x_{i, \alpha \beta}x_{j}) \triangleright (\fa_- \shuffle  \fb_-)\\
    =& \ \left(x_{a + b, \alpha \beta} + \sum \limits_{i+j = a + b} \Delta^j_{a,b} x_{i, \alpha \beta}x_{j}\right) \triangleright (\fa_- \shuffle  \fb_-)\\
    =& \ (x_{a, \alpha} \diamond x_{b, \beta}) \triangleright (\fa_- \shuffle  \fb_-).
\end{align*}
This proves part $(1)$. Part $(2)$ is straightforward from the commutativity of the shuffle product (see Proposition \ref{prop: commutative}). This completes the proof.
\end{proof}

We now ready to prove the main result of this subsection.

\begin{proposition} \label{prop: associative}
  The diamond product and the shuffle product on $\mathfrak{D}$ are associative.
\end{proposition}

\begin{proof}
Let $\fa, \fb, \fc \in \langle \Gamma \rangle$ be arbitrary words. It suffices to show that
\begin{equation*} \label{eq: associative}
    (\fa \diamond \fb) \diamond \fc = \fa \diamond (\fb \diamond \fc) \quad \text{and} \quad (\fa \shuffle \fb) \shuffle \fc = \fa \shuffle (\fb \shuffle\fc) .
\end{equation*}
We proceed the proof by induction on $\text{depth}(\fa) + \text{depth}(\fb) + \text{depth}(\fc)$. If one of $\fa, \fb$ or $\fc$ is the empty word, then \eqref{eq: associative} holds trivially. We assume that \eqref{eq: associative} holds when $\text{depth}(\fa) + \text{depth}(\fb) + \text{depth}(\fc) < n$ with $n \in \mathbb{N}$ and $n \geq 3$. We need to show that \eqref{eq: associative} holds when $\text{depth}(\fa) + \text{depth}(\fb) + \text{depth}(\fc) = n$.

We first show that $(\fa \diamond \fb) \diamond \fc = \fa \diamond (\fb \diamond \fc)$. Assume that
$\fa = x_{a,\alpha} \fa_-,\fb = x_{b,\beta} \fb_- = \fc = x_{c,\gamma} \fc_-$. From Lemma \ref{lem: triangleright formulas}, we have
\begin{align*}
    (\fa &\diamond \fb) \diamond \fc \\
    =& \ \left(x_{a + b, \alpha \beta}(\fa_- \shuffle \fb_-) + \sum\limits_{i+j = a + b} \Delta^j_{a,b} x_{i, \alpha \beta}(x_{j} \shuffle (\fa_- \shuffle \fb_-))\right) \diamond \fc\\
    =& \ x_{a + b, \alpha \beta}(\fa_- \shuffle \fb_-) \diamond \fc + \sum\limits_{i+j = a + b} \Delta^j_{a,b} x_{i, \alpha \beta}(x_{j} \shuffle (\fa_- \shuffle \fb_-))  \diamond \fc\\
    =& \ (x_{a + b, \alpha \beta} \diamond x_{c, \gamma}) \triangleright ((\fa_-\shuffle \fb_-) \shuffle \fc_-)\\
    &+ \sum\limits_{i+j = a + b} \Delta^j_{a,b} (x_{i, \alpha \beta} \diamond x_{c, \gamma}) \triangleright ((x_{j} \shuffle (\fa_- \shuffle \fb_-)) \shuffle \fc_-).
\end{align*}
For all $i,j \in \mathbb{N}$ with $i+j = a+b$, it follows from the induction hypothesis that
\begin{align*}
    & (x_{i, \alpha \beta} \diamond x_{c, \gamma}) \triangleright ((x_{j} \shuffle (\fa_- \shuffle \fb_-)) \shuffle \fc_-)\\
    =& \ (x_{i, \alpha \beta} \diamond x_{c, \gamma}) \triangleright (x_{j} \shuffle ((\fa_- \shuffle \fb_-) \shuffle \fc_-))\\
    =& \ \left(x_{i + c,\alpha \beta\gamma} + \sum\limits_{i_1+j_1 = i + c} \Delta^{j_1}_{i,c} x_{i_1,\alpha \beta\gamma}x_{j_1}\right) \triangleright (x_{j} \shuffle ((\fa_- \shuffle \fb_-) \shuffle \fc_-))\\
    =& \ x_{i + c,\alpha \beta\gamma}\triangleright (x_{j} \shuffle ((\fa_- \shuffle \fb_-) \shuffle \fc_-)) \\
    &+ \sum\limits_{i_1+j_1 = i + c} \Delta^{j_1}_{i,c} x_{i_1,\alpha \beta\gamma}x_{j_1} \triangleright (x_{j} \shuffle ((\fa_- \shuffle \fb_-) \shuffle \fc_-))\\
    =& \ x_{i + c,\alpha \beta\gamma} (x_{j} \shuffle ((\fa_- \shuffle \fb_-) \shuffle \fc_-)) \\
    &+ \sum\limits_{i_1+j_1 = i + c} \Delta^{j_1}_{i,c} x_{i_1,\alpha \beta\gamma} ((x_{j_1} \shuffle x_{j} ) \shuffle ((\fa_- \shuffle \fb_-) \shuffle \fc_-)) \\
    =& \ (x_{i + c,\alpha \beta\gamma} x_{j}) \triangleright ((\fa_- \shuffle \fb_-) \shuffle \fc_-) \\
    &+ \sum\limits_{i_1+j_1 = i + c} \Delta^{j_1}_{i,c} x_{i_1,\alpha \beta\gamma}(x_{j_1} \shuffle x_{j}) \triangleright  ((\fa_- \shuffle \fb_-) \shuffle \fc_-)\\
    =& \left( x_{i + c,\alpha \beta\gamma} x_{j} + \sum\limits_{i_1+j_1 = i + c} \Delta^{j_1}_{i,c} x_{i_1,\alpha \beta\gamma}(x_{j_1} \shuffle x_{j})\right) \triangleright  ((\fa_- \shuffle \fb_-) \shuffle \fc_-)\\
    =& \ ((x_{i ,\alpha \beta}x_{j}) \diamond x_{c,\gamma}) \triangleright ((\fa_- \shuffle \fb_-) \shuffle \fc_-).
\end{align*}
Thus
\begin{align*}
    (\fa &\diamond \fb) \diamond \fc\\
    =& \ (x_{a + b, \alpha \beta} \diamond x_{c, \gamma}) \triangleright ((\fa_-\shuffle \fb_-) \shuffle \fc_-) \\
    &+ \sum\limits_{i+j = a + b} \Delta^j_{a,b} ((x_{i ,\alpha \beta}x_{j}) \diamond x_{c,\gamma}) \triangleright ((\fa_- \shuffle \fb_-) \shuffle \fc_-)\\
    =& \left((x_{a + b, \alpha \beta} \diamond x_{c, \gamma}) + \sum\limits_{i+j = a + b} \Delta^j_{a,b} (x_{i ,\alpha \beta}x_{j} \diamond x_{c,\gamma}) \right)\triangleright ((\fa_- \shuffle \fb_-) \shuffle \fc_-)\\
    =& \ ((x_{a, \alpha} \diamond x_{b, \beta}) \diamond x_{c, \gamma}) \triangleright ((\fa_- \shuffle \fb_-) \shuffle \fc_-).
\end{align*}
On the other hand, from Proposition \ref{prop: commutative} and the above arguments, we deduce that
\begin{align*}
    \fa \diamond (\fb \diamond \fc) &= (\fb \diamond \fc) \diamond \fa  \\
    &=  ((x_{b, \beta} \diamond x_{c, \gamma}) \diamond x_{a, \alpha})\triangleright ((\fb_- \shuffle \fc_-) \shuffle \fa_-) \\
    &= (x_{a, \alpha} \diamond (x_{b, \beta} \diamond x_{c, \gamma})) \triangleright (\fa_- \shuffle  (\fb_- \shuffle \fc_-)).
\end{align*}
It follows from Proposition \ref{prop: associativiy of letters} that $(x_{a,\alpha} \diamond x_{b,\beta}) \diamond x_{c,\gamma} = x_{a,\alpha} \diamond (x_{b,\beta} \diamond x_{c,\gamma})$. Moreover, it follows from the induction hypothesis that $(\fa_- \shuffle \fb_-) \shuffle \fc_- = \fa_- \shuffle  (\fb_- \shuffle \fc_-)$. We thus conclude that $(\fa \diamond \fb) \diamond \fc = \fa \diamond (\fb \diamond \fc)$.

Next, we show that $(\fa \shuffle \fb) \shuffle \fc = \fa \shuffle (\fb \shuffle \fc)$. From Lemma \ref{lem: triangleright formulas}, we have
\begin{align*}
    (\fa \shuffle \fb) \shuffle \fc
    =& \ (\fa \triangleright \fb + \fb \triangleright \fa + \fa \diamond \fb) \shuffle \fc\\
    =& \ (\fa \triangleright \fb) \shuffle \fc  + (\fb \triangleright \fa) \shuffle \fc  + (\fa \diamond \fb) \shuffle \fc\\
    =& \ ((\fa \triangleright \fb) \triangleright \fc + \fc \triangleright (\fa \triangleright \fb) + (\fa \triangleright \fb) \diamond \fc)\\
    &+  ((\fb \triangleright \fa) \triangleright \fc + \fc \triangleright (\fb \triangleright \fa) + (\fb \triangleright \fa) \diamond \fc)\\
    &+  ((\fa \diamond \fb) \triangleright \fc + \fc \triangleright (\fa \diamond \fb) + (\fa \diamond \fb) \diamond \fc)\\
    =& \ (\fa \triangleright \fb) \triangleright \fc + (\fa \triangleright \fb) \diamond \fc + (\fb \triangleright \fa) \triangleright \fc + (\fb \triangleright \fa) \diamond \fc + (\fa \diamond \fb) \triangleright \fc\\
    &+ (\fc \triangleright (\fa \triangleright \fb) + \fc \triangleright (\fb \triangleright \fa) + \fc \triangleright (\fa \diamond \fb)) + (\fa \diamond \fb) \diamond \fc\\
    =& \ (\fa \triangleright \fb) \triangleright \fc + (\fa \triangleright \fb) \diamond \fc + (\fb \triangleright \fa) \triangleright \fc + (\fb \triangleright \fa) \diamond \fc + (\fa \diamond \fb) \triangleright \fc\\
    &+ \fc \triangleright (\fa \shuffle \fb)  + (\fa \diamond \fb) \diamond \fc
\end{align*}
and
\begin{align*}
    \fa \shuffle (\fb \shuffle \fc)
    =& \ \fa \shuffle (\fb \triangleright \fc + \fc \triangleright \fb + \fb \diamond \fc)\\
    =& \ \fa \shuffle (\fb \triangleright \fc) + \fa \shuffle (\fc \triangleright \fb) + \fa \shuffle (\fb \diamond \fc)\\
    =& \ (\fa \triangleright (\fb \triangleright \fc) + (\fb \triangleright \fc) \triangleright \fa + \fa \diamond (\fb \triangleright \fc))\\
    &+ (\fa \triangleright (\fc \triangleright \fb) + (\fc \triangleright \fb) \triangleright \fa + \fa \diamond (\fc \triangleright \fb)) \\
    &+ (\fa \triangleright (\fb \diamond \fc) + (\fb \diamond \fc) \triangleright \fa + \fa \diamond (\fb \diamond \fc))\\
    =& \ (\fb \triangleright \fc) \triangleright \fa + \fa \diamond (\fb \triangleright \fc) + (\fc \triangleright \fb) \triangleright \fa + \fa \diamond (\fc \triangleright \fb) + (\fb \diamond \fc) \triangleright \fa\\
    &+ (\fa \triangleright (\fb \triangleright \fc) + \fa \triangleright (\fc \triangleright \fb) + \fa \triangleright (\fb \diamond \fc)) + \fa \diamond (\fb \diamond \fc)\\
    =& \ (\fb \triangleright \fc) \triangleright \fa + \fa \diamond (\fb \triangleright \fc) + (\fc \triangleright \fb) \triangleright \fa + \fa \diamond (\fc \triangleright \fb) + (\fb \diamond \fc) \triangleright \fa\\
    &+ \fa \triangleright (\fb \shuffle \fc) + \fa \diamond (\fb \diamond \fc).
\end{align*}

We now compare the above expansions. We have showed that $(\fa \diamond \fb) \diamond \fc = \fa \diamond (\fb \diamond \fc)$. On the other hand, we have
\begin{equation*}
    \fc \triangleright (\fa \shuffle \fb) = x_{c,\gamma}(\fc_- \shuffle(\fa \shuffle \fb))
    \quad \text{and} \quad
    (\fc \triangleright \fb) \triangleright \fa = x_{c,\gamma}(\fc_- \shuffle \fb) \triangleright \fa = x_{c,\gamma}((\fc_- \shuffle \fb) \shuffle \fa).
\end{equation*}
From the induction hypothesis and commutativity of shuffle product, one deduces that $\fc \triangleright (\fa \shuffle \fb) = (\fc \triangleright \fb) \triangleright \fa$. Similarly, one deduces that $(\fa \triangleright \fb) \triangleright \fc= \fa \triangleright (\fb \shuffle \fc)$.

We have
\begin{align*}
    &(\fa \triangleright \fb) \diamond \fc
    = x_{a,\alpha}(\fa_- \shuffle \fb) \diamond \fc
    = (x_{a,\alpha} \diamond x_{c,\gamma}) \triangleright ((\fa_- \shuffle \fb) \shuffle \fc_-), \\
    & \fa \diamond (\fc \triangleright \fb) = \fa \diamond x_{c,\gamma}(\fc_- \shuffle \fb)
    =  (x_{a,\alpha} \diamond x_{c,\gamma}) \triangleright (\fa_- \shuffle (\fc_- \shuffle \fb)).
\end{align*}

From the induction hypothesis and commutativity of shuffle product, one deduces that $(\fa \triangleright \fb) \diamond \fc = \fa \diamond (\fc \triangleright \fb)$.

We have
\begin{equation*}
\begin{split}
    (\fb \triangleright \fa) \triangleright \fc =& \ x_{b,\beta}(\fb_- \shuffle \fa) \triangleright \fc\\
    =& \  x_{b,\beta}((\fb_- \shuffle \fa) \shuffle \fc)\\
\end{split}
     \quad \text{and} \quad
\begin{split}
    (\fb \triangleright \fc) \triangleright \fa =& \ x_{b,\beta}(\fb_- \shuffle \fc) \triangleright \fa\\
    =& \ x_{b,\beta}((\fb_- \shuffle \fc) \shuffle \fa).
\end{split}
\end{equation*}
From the induction hypothesis and commutativity of shuffle product, one deduces that $ (\fb \triangleright \fa) \triangleright \fc = (\fb \triangleright \fc) \triangleright \fa$.

It follows from the induction hypothesis and that
\begin{align*}
(\fa \diamond \fb) \triangleright \fc
    =& \ \left(x_{a + b, \alpha \beta}(\fa_- \shuffle \fb_-) + \sum\limits_{i+j = a + b} \Delta^j_{a,b} x_{i, \alpha \beta}(x_{j} \shuffle (\fa_- \shuffle \fb_-))\right) \triangleright \fc\\
    =& \ x_{a + b, \alpha \beta}(\fa_- \shuffle \fb_-)\triangleright \fc + \sum\limits_{i+j = a + b} \Delta^j_{a,b} x_{i, \alpha \beta}(x_{j} \shuffle (\fa_- \shuffle \fb_-)) \triangleright \fc\\
    =& \ x_{a + b, \alpha \beta}((\fa_- \shuffle \fb_-) \shuffle \fc) + \sum \limits_{i+j = a + b} \Delta^j_{a,b} x_{i, \alpha \beta} ((x_{j} \shuffle (\fa_- \shuffle \fb_-)) \shuffle\fc) \\
    =& \ x_{a + b, \alpha \beta}((\fa_- \shuffle \fb_-) \shuffle \fc) + \sum \limits_{i+j = a + b} \Delta^j_{a,b} x_{i, \alpha \beta} (x_{j} \shuffle ((\fa_- \shuffle \fb_-) \shuffle\fc)) \\
    =& \ x_{a + b, \alpha \beta} \triangleright ((\fa_- \shuffle \fb_-) \shuffle \fc) + \sum \limits_{i+j = a + b} \Delta^j_{a,b} x_{i, \alpha \beta}x_{j} \triangleright ((\fa_- \shuffle \fb_-) \shuffle\fc)\\
    =& \ \left(x_{a + b, \alpha \beta} + \sum \limits_{i+j = a + b} \Delta^j_{a,b} x_{i, \alpha \beta}x_{j} \right) \triangleright ((\fa_- \shuffle \fb_-) \shuffle\fc)\\
    =& \ (x_{a,\alpha} \diamond x_{b,\beta}) \triangleright ((\fa_- \shuffle \fb_-) \shuffle\fc)\\
\end{align*}
and
\begin{align*}
    \fa \diamond (\fb \triangleright \fc) = \fa \diamond x_{b,\beta}(\fb_- \shuffle \fc) = (x_{a,\alpha} \diamond x_{b,\beta}) \triangleright (\fa_- \shuffle (\fb_- \shuffle \fc)).
\end{align*}
From the induction hypothesis, one deduces that $(\fa \diamond \fb) \triangleright \fc = \fa \diamond (\fb \triangleright \fc)$. Similarly, one deduces that $(\fb \triangleright \fa) \diamond \fc = (\fb \diamond \fc) \triangleright \fa$.

From the above arguments, we conclude that $(\fa \shuffle \fb) \shuffle \fc = \fa \shuffle (\fb \shuffle \fc)$. This completes the proof.
\end{proof}

As for $\mathfrak{C}$, we define the unit $u \colon \mathbb{F}_q\to \mathfrak{D}$, sending $1 \in \mathbb{F}_q$ to an empty word $1 \in \mathfrak{D}$. Then,
as a direct consequence of Proposition \ref{prop: commutative} and Proposition \ref{prop: associative}, we obtain the following result.

\begin{theorem} \label{thm: commutative algebra}
The spaces $(\mathfrak{D}, \diamond)$ and $(\mathfrak{D}, \shuffle )$ are commutative $\Fq$-algebras.
\end{theorem}

The following proposition summarizes several properties of different products $\triangleright$, $\diamond$ and $\shuffle $ that are shown above and will be useful in the sequel.

\begin{proposition} \label{prop: properties on triangle diamond sha}
For all nonempty words $\mathfrak{a}, \mathfrak{b} , \mathfrak{c}\in \langle \Gamma \rangle$, we have
\begin{enumerate}[\normalfont (1)]
    \item $(\mathfrak{a} \triangleright \mathfrak{b}) \triangleright \mathfrak{c} = \mathfrak{a} \triangleright (\mathfrak{b} \shuffle \mathfrak{c})$,
    \item $(\mathfrak{a} \triangleright \mathfrak{b}) \diamond \mathfrak{c} = \mathfrak{a} \diamond (\mathfrak{c} \triangleright \mathfrak{b}) = (\mathfrak{a} \diamond \mathfrak{c}) \triangleright \mathfrak{b}$.
\end{enumerate}
\end{proposition}


\subsection{Shuffle map for AMZV's in positive characteristic} ${}$\par

For each $d\in \mathbb{Z}$, we define the following two $\mathbb{F}_q$-linear maps
\[S_{<d} \colon \mathfrak{D}\to K_\infty \quad
\text{and}
\quad
\zeta_A\colon \mathfrak{D}\to K_\infty,
\] both of which map the empty word $1 \in \mathfrak{D}$ to $1 \in K_\infty$, and the word $x_{s_1, \varepsilon_1} \dots x_{s_r, \varepsilon_r}$ to
$ S_{<d} \begin{pmatrix}
 \varepsilon_1 & \dotsb & \varepsilon_n \\
s_1 & \dotsb & s_n \end{pmatrix}$ and $\zeta_A\begin{pmatrix}
 \varepsilon_1 & \dotsb & \varepsilon_n \\
s_1 & \dotsb & s_n \end{pmatrix}$ respectively, identifying presentations
\begin{align*}
S_{<d} (x_{s_1, \varepsilon_1}\dots x_{s_r, \varepsilon_r})&=S_{<d}\begin{pmatrix}
 \varepsilon_1 & \dotsb & \varepsilon_n \\
s_1 & \dotsb & s_n \end{pmatrix},\\
\zeta_A (x_{s_1, \varepsilon_1}\dots x_{s_r, \varepsilon_r})& =\zeta_A\begin{pmatrix}
 \varepsilon_1 & \dotsb & \varepsilon_n \\
s_1 & \dotsb & s_n \end{pmatrix}.
\end{align*}

\begin{theorem} \label{thm: AZ shuffle map}
    For $\fa, \fb \in \mathfrak{D}$
    and $d \in \mathbb{Z}$ we have
    \begin{align*}
    S_{<d}(\fa \sha \fb) &= S_{<d}(\fa)S_{<d}(\fb).\\
    \zeta_{A}(\fa \sha \fb) &= \zeta_{A}(\fa)\zeta_{A}(\fb).
    \end{align*}
\end{theorem}

\begin{proof}
The $\shuffle$ product is defined to satisfy this equation. We can check that the $\shuffle$ follows the inductive steps for product of power sums and $\zeta_A$ as follows.

    When any of $\fa$ and $\fb$ is the empty word, the theorem holds trivially. Let $\fa = x_{a, \alpha} \fa_-$ and $\fb = x_{b, \beta} \fb_-$ be nonempty words. We proceed the proof by induction on $\depth(\fa)+ \depth(\fb)$. Following the proof of \cite[Lemma~2.5]{Har21} and applying Lemma~\ref{lemma: VMPLdepth1},
    \begin{align*}
         S_{<d}(\fa) S_{<d}(\fb)
        = &
         \sum_{m <d }S_{m}(x_{a, \alpha}) S_{< m}(\fa_-)S_{< m}(\fb) + \sum_{m <d }S_{m}(x_{b, \beta}) S_{< m}(\fa)S_{< m}(\fb_-) \\
         & + \sum_{m<d} S_{m}(x_{a, \alpha}) S_{m}(x_{b, \beta})
         S_{<m}(\fa_-) S_{<m}(\fb_-)\\
        =&
        \sum_{m <d }S_{m}(x_{a, \alpha}) S_{< m}(\fa_-)S_{< m}(\fb) + \sum_{m <d }S_{m}(x_{b, \beta}) S_{< m}(\fa)S_{< m}(\fb_-) \\
        & + \sum_{m<d}
        S_{m}(x_{a+b, \alpha\beta})
        S_{<m}(\fa_-) S_{<m}(\fb_-)\\
        & + \sum_{m<d} \Delta_{a, b}^j
        \sum_{i+j = a+b}  S_{m}(x_{i, \alpha\beta} x_{j})
        S_{<m}(\fa_-) S_{<m}(\fb_-).
    \end{align*}
    This is parallel to the recursive definition of $\sha$. Proof is the same for $\zeta_A$.
\end{proof}


\section{Shuffle Hopf algebra of AMZV's} \label{sec: shuffle Hopf algebra}

In \S \ref{sec: shuffle algebra} we have defined the shuffle algebra structure of AMZV's. This section is devoted to the construction of the Hopf algebra structure of AMZV's. In \S \ref{sec: coproduct} we present the coproduct for the word algebra $\frak D$. In \S \ref{sec: compatibility} we prove the compatibility of the coproduct in Theorem \ref{thm: main statement on compatibility}. Unlike the MZV's case, the proof of this theorem uses the induction on depth instead of that on weight. Then in \S \ref{sec: coassociativity} we prove the coassociativity of the coproduct in Theorem \ref{thm: coassociativity}. Putting all together we get the main result of this paper (see Theorem \ref{thm: AMZV Hopf algebra}).

\subsection{Coproduct} \label{sec: coproduct} ${}$\par

This section is inspired by our previous work \cite[\S 8]{IKLNDP23}. We will introduce the coproduct
	\[ \Delta: \frak D \to \frak D \otimes \frak D. \]

First, we define the coproduct of depth one. We put
\begin{align*}
\Delta(1):=1 \otimes 1, \quad\quad
\Delta(x_{1,\varepsilon}):=1 \otimes x_{1,\varepsilon} + x_{1,\varepsilon} \otimes 1
\end{align*}
for all $\varepsilon \in \Fq^*$.

Further, we define the coproduct of depth one as follows.
\begin{definition} \label{defn: depth one coproduct}
For all $n \in \N$ and $\varepsilon \in \Fq^*$, we set
\begin{equation*}
    \Delta(x_{n,\varepsilon}) = 1 \otimes x_{n,\varepsilon} + \sum \limits_{\substack{r \in \N, \fa \in \langle \Sigma \rangle\\ r + w(\fa) = n}} \binom{r + \depth(\fa) - 2}{\depth(\fa)} x_{r,\varepsilon} \otimes [\fa].
\end{equation*}
Here we recall that $[\fa]$ is given as in \eqref{eq: bracket operator} and we identify $\fa \in \Sigma$ with $\varphi_1(\fa) \in \Gamma$.
\end{definition}

Next, we define the coproduct for words of depth greater than one. Let $n \in \bN$ with $n>1$ and we suppose that we have defined $\Delta(\fu)$ for all words $\fu$ of $\depth(\fu)<n$. We now give a formula for $\Delta(\fu)$ for all words~$\fu$ with $\depth(\fu)=n$. For such a word~$\fu$, we put $\fu=x_{u,\varepsilon} \fv$ with $\depth(\fv)=n-1<n$. Since $x_{u,\varepsilon}$ and $\fv$ are both of depth strictly less than $n$, we have already defined
\begin{align*}
\Delta(x_{u,\varepsilon})=1 \otimes x_{u,\varepsilon} + \sum a_{u,\varepsilon} \otimes b_{u,\varepsilon}, \text{ and }\quad \Delta(\fv)= \sum a_\fv \otimes b_\fv,
\end{align*}
with $a_{u,\varepsilon} \in \Gamma$ and $b_{u,\varepsilon}, a_\fv, b_\fv \in \langle \Gamma \rangle$. Then we set
\begin{align*}
\Delta(\fu) &:=1 \otimes \fu + \sum (a_{u,\varepsilon} \triangleright a_\fv) \otimes (b_{u,\varepsilon} \shuffle  b_\fv) \\
& =1 \otimes \fu + \sum (a_{u,\varepsilon} a_\fv) \otimes (b_{u,\varepsilon} \shuffle  b_\fv).
\end{align*}
Here the last equality holds as $a_{u,\varepsilon} \in \Gamma$.

We end this section by proving the following result:

\begin{proposition} \label{prop: Delta varphi}
For any nonempty word $\fu  \in \langle \Gamma \rangle$, we have
\begin{equation*}
\Delta \left( \varphi_\varepsilon (\fu)\right) = ( \varphi_\varepsilon \otimes \Id ) \left( \Delta ( \fu)\right) + ( \Id \otimes \varphi_\varepsilon - \varphi_\varepsilon \otimes \Id) (1\otimes \fu).
\end{equation*}
\end{proposition}

\begin{proof}
Note that for $\fa = x_a \fa_- \in \langle \Gamma \rangle$, $\fb \in \langle \Gamma\rangle$ and $\varepsilon \in \Fq^*$,
	\[ \varphi_\varepsilon(\fa) \triangleright \fb = x_{a, \varepsilon} ( \fa_- \shuffle \fb) = \varphi_\varepsilon( \fa \triangleright \fb). \]
We mention that this does not hold when $\fa =1$.

We first consider the case where $\depth(\fu)=1$, i.e., $\fu$ belongs to $\Gamma$. We put $\fu=x_{n,\alpha}$ with $n \in \mathbb N$ and $\alpha \in \Fq^*$. Thus $\varphi_\varepsilon(\fu)=x_{n,\varepsilon \alpha}=\varphi_{\varepsilon \alpha}(x_n)$.

Recall that (see Proposition \ref{prop: formula delta xn})
\begin{equation*}
    \Delta(x_{n}) = 1 \otimes x_{n} + \sum \limits_{\substack{r \in \N, \fa \in \langle \Sigma \rangle\\ r + w(\fa) = n}} \binom{r + \depth(\fa) - 2}{\depth(\fa)} x_{r} \otimes [\fa].
\end{equation*}
Therefore, Definition~\ref{defn: depth one coproduct} can be written as
\begin{align*}
\Delta \left(\fu\right)&  = \Delta \left( \varphi_\alpha (x_n)\right) \\
&  =  \left( \varphi_\alpha \otimes \Id \right) \left( \Delta(x_n) - 1\otimes x_n\right) + \left( \Id  \otimes \varphi_\alpha \right)(1 \otimes x_{n}) \\
& =  \left( \varphi_\alpha \otimes \Id \right) \left( \Delta(x_n)\right)  + (\Id \otimes \varphi_\alpha - \varphi_\alpha \otimes \Id )\left(1\otimes x_{n}\right).
\end{align*}
Similarly,
\begin{align*}
\Delta \left(\varphi_\varepsilon(\fu)\right)&  = \Delta \left(\varphi_{\varepsilon \alpha}(x_n)\right) \\
&  =  \left( \varphi_{\varepsilon \alpha} \otimes \Id \right) \left( \Delta(x_n) - 1\otimes x_n\right) + \left( \Id  \otimes \varphi_{\varepsilon \alpha} \right)(1 \otimes x_{n}) \\
& =  \left( \varphi_{\varepsilon \alpha} \otimes \Id \right) \left( \Delta(x_n)\right)  + (\Id  \otimes \varphi_{\varepsilon \alpha} - \varphi_{\varepsilon \alpha} \otimes \Id)\left(1\otimes x_{n}\right).
\end{align*}
By Lemma \ref{lemma: associativity horizontal maps}, we conclude that for $\fu=x_{n,\alpha}$,
\begin{equation*}
\Delta \left( \varphi_\varepsilon (\fu)\right) = ( \varphi_\varepsilon \otimes \Id ) \left( \Delta ( \fu)\right) + ( \Id  \otimes \varphi_\varepsilon - \varphi_\varepsilon \otimes \Id ) (1\otimes \fu).
\end{equation*}

Next, suppose that $\fu = x_{u,\alpha} \fv \in \langle \Gamma \rangle$ with $\depth(\fu) \ge 2$. Let $\Delta(x_u) = 1 \otimes x_u + \sum a_u \otimes b_u$, and $\Delta(\fv) = \sum a_\fv \otimes b_\fv$. We have $a_u \in \Sigma$, $b_u \in \langle \Sigma\rangle$ with $a_u \ne 1$ by Lemma~\ref{lem: au non empty}, and $a_\fv, b_\fv \in \langle \Gamma \rangle$. We have seen that
	\[ \Delta(x_{u,\alpha}) = 1 \otimes x_{u,\alpha} + \sum a_{u,\alpha} \otimes b_u. \]
By definition,
\begin{align*}
    \Delta(\fu) &=  1\otimes \fu + \sum (a_{u,\alpha} \triangleright a_\fv ) \otimes
    (b_u \shuffle b_\fv) \\
    &=  1\otimes \fu + \sum (a_{u,\alpha} a_\fv ) \otimes
    (b_u \shuffle b_\fv).
\end{align*}
The last equality follows from the fact that $\depth(a_{u,\alpha})=1$.

Since $\varphi_\varepsilon (\fu) = x_{u, \varepsilon \alpha} \fv$ and $\Delta\left( \varphi_\varepsilon(x_{u,\alpha})\right) =
1\otimes \varphi_\varepsilon (x_{u,\alpha}) + \sum \varphi_\varepsilon (a_{u,\alpha}) \otimes b_u
$,
\begin{align*}
\Delta\left(\varphi_\varepsilon(\fu)\right) &= \Delta\left(\varphi_\varepsilon(x_{u,\alpha}) \fv \right) \\
& =1 \otimes \varphi_\varepsilon(\fu) + \sum \left(\varphi_\varepsilon(a_{u,\alpha}) a_\fv\right) \otimes (b_{u} \shuffle  b_\fv)\\
& = 1 \otimes \varphi_\varepsilon(\fu) + \sum\varphi_\varepsilon \left( a_{u,\alpha} a_\fv\right) \otimes (b_{u} \shuffle  b_\fv).
\end{align*}
Thus we get
\begin{equation*}
\Delta \left( \varphi_\varepsilon (\fu)\right) = ( \varphi_\varepsilon \otimes \Id ) \left( \Delta ( \fu)\right) + ( \Id  \otimes \varphi_\varepsilon - \varphi_\varepsilon \otimes\Id ) (1\otimes \fu).
\end{equation*}
\end{proof}

\subsection{Compatibility} \label{sec: compatibility}

\subsubsection{Compatibility for letters} ${}$\par

We want to prove the following result (for MZV's, see \cite[Proposition 7.5]{IKLNDP23}):

\begin{proposition} \label{prop: compatibility step 3}
Let $u,v \in \bN$ and $\alpha, \beta \in \Fq^*$.
Then
	\[ \Delta(x_{u,\alpha} \shuffle x_{v,\beta})=\Delta(x_{u,\alpha})\shuffle \Delta(x_{v,\beta}). \]
\end{proposition}
\begin{proof}
    Let \[ \Delta(x_u) = 1\otimes x_u + \sum a_u \otimes b_u, \quad \Delta(x_v) = 1\otimes x_v + \sum a_v \otimes b_v.\]
    By Proposition~\ref{prop: Delta varphi}, we have
    \[ \Delta(x_{u, \alpha}) = 1\otimes x_{u, \alpha} + \sum a_{u, \alpha} \otimes b_u, \quad \Delta(x_{v, \beta}) = 1\otimes x_{v, \beta} + \sum a_{v, \beta} \otimes b_v,\]
    where $a_{u, \alpha} = \varphi_\alpha(a_u)$ and $a_{v, \beta} = \varphi_\beta(a_v)$. Here we recall that $a_u$ and $a_v$ are of depth~$1$.

    Since $x_{u, \alpha} \shuffle x_{v, \beta}= x_{u, \alpha} x_{v, \beta} + x_{v, \beta}x_{u, \alpha} + \varphi_{\alpha \beta} (x_u \diamond  x_v)$ (by Lemma~\ref{lem: distributive}),
    \begin{align*}
        \Delta\left( x_{u, \alpha} \shuffle x_{v, \beta}\right)  =&
        1 \otimes x_{u, \alpha} x_{v, \beta} + \sum a_{u, \alpha} \otimes (b_u \shuffle x_{v, \beta}) +
        \sum (a_{u, \alpha } a_{v, \beta}) \otimes (b_u \shuffle b_v) \\
        & + 1 \otimes x_{v, \beta} x_{u, \alpha} + \sum a_{v, \beta} \otimes (b_v \shuffle x_{u, \alpha}) +
        \sum (a_{v, \beta} a_{u, \alpha}) \otimes (b_v \shuffle b_u) \\
        & + \Delta(\varphi_{\alpha \beta} (x_{u}\diamond x_v)).
    \end{align*}
    On the other hand,
    \begin{align*}
        \Delta( x_{u, \alpha}) \shuffle \Delta(x_{v, \beta}) =&
        1 \otimes (x_{u, \alpha} \shuffle x_{v, \beta}) +
        \sum a_{u, \alpha} \otimes (b_u \shuffle x_{v, \beta})\\
        & + \sum a_{v, \beta} \otimes (b_v \shuffle x_{u, \alpha})
        + \sum (a_{u, \alpha} \shuffle a_{v, \beta}) \otimes
        (b_{u} \shuffle b_v).
    \end{align*}
    By expanding $\shuffle$ in each tensorand and canceling, the proposition holds when
    \begin{align*} \Delta(\varphi_{\alpha\beta}(x_u \diamond x_v))
    & =
    1 \otimes \varphi_{\alpha \beta} (x_u \diamond x_v) +
    \sum (a_{u, \alpha} \diamond a_{v, \beta})\otimes (b_u \shuffle b_v) \\
    & = 1 \otimes \varphi_{\alpha \beta} (x_u \diamond x_v) +
    ( \varphi_{\alpha \beta} \otimes \Id)
    \left( \sum (a_{u} \diamond a_{v})\otimes (b_u \shuffle b_v) \right).
    \end{align*}
    By Proposition \ref{prop: Delta varphi}, the proposition is reduced into
    \begin{align*}
    \Delta( x_u \diamond x_v)  =   1\otimes (x_u \diamond x_v)
    +
     \sum (a_{u} \diamond a_{v})\otimes (b_u \shuffle b_v),
    \end{align*}
    which is a special case of Lemma~\ref{lem: delta diamond}.
\end{proof}

\subsubsection{Compatibility for words of arbitrary depth} ${}$\par

\begin{lemma} \label{lem: sha and depth}
    Let $\fu, \fv \in \langle \Gamma \rangle$. Then $\fu \sha \fv$ can be presented as the sum of words whose depths are all less than or equal to $\depth(\fu)+\depth(\fv)$.
\end{lemma}

\begin{proof}

We proceed with the induction on the total depth $\depth(\fu)+\depth(\fv)$.

The lemma is trivial when $\fu$ or $\fv$ is trivial. If $\fu = x_{u, \alpha}$ and $\fv = x_{v, \beta}$, then the lemma is also immediate from the expression
    \[ \fu \sha \fv = x_{u+v, \alpha\beta}  + x_{u, \alpha}x_{v,\beta} + x_{v, \beta}x_{u,\alpha} + \sum_{i+j=u+v} \Delta_{u, v}^j x_{i, \alpha\beta} x_j. \]
So we assume that $\depth(\fu), \depth(\fv)\ge1$ with at least one of $\fu$ and $\fv$ has depth~$>1.$

    Now assume the lemma holds for all $\fa,\fb \in \langle \Gamma \rangle$ such that $\depth(\fa)+\depth(\fb)<\depth(\fu)+\depth(\fv)$. Let $\fu = x_{u, \alpha}\fu_-$, $\fv = x_{v, \beta}\fv_-$. Then
    \begin{align*}
        x_{u, \alpha}\fu_- \sha x_{v, \beta}\fv_- & = x_{u+v, \alpha \beta}(\fu_- \sha \fv_-) + x_{u, \alpha}(\fu_- \sha \fv ) + x_{v, \beta}(\fu \sha \fv_-) \\
        & + \sum_{i+j=u+v}\Delta_{u, v}^j x_{i, \alpha\beta} (x_j \sha \fu_- \sha \fv_-)
    \end{align*}
    From the induction hypothesis, the first three terms can be presented as a sum of words with depths less than or equal to $\depth(\fu) + \depth(\fv)$. For each $j$, $x_j \sha \fu_- = \fu_- \sha x_j$ can be presented as a sum of words whose depths are less than or equal to $\depth(\fu)$ from the induction hypothesis. Therefore, $x_i(x_j \sha \fu_- \sha \fv_-)$ can be presented as sum of words with depth $\depth(\fu)+ \depth(\fv)$. This proves the lemma.
\end{proof}

The proof of the compatibility of $\Delta$ and $\sha$ for $\langle \Gamma \rangle$ is parallel to that for $\langle \Sigma \rangle$.
\begin{theorem} \label{thm: main statement on compatibility}
    Let $\fu, \fv\in \langle \Gamma\rangle$. Then $\Delta$ and $\sha$ are compatible for $\fu, \fv$, that means
    \[ \Delta(\fu \sha \fv) =\Delta(\fu)\sha \Delta(\fv). \]
\end{theorem}

 Note that the theorem is immediate when $\fu$ or $\fv$ is trivial. Also, $\depth(\fu)=\depth(\fv)=1$ case is shown in Proposition~\ref{prop: compatibility step 3}. Before proving this theorem, we need some preparatory results.

 \begin{lemma} \label{lemma: compatibility first character reduce}
Let $\fu = x_{u, \alpha} \fu_-$, $\fv= x_{v, \beta} \fv_-$ be two words in $\langle \Gamma \rangle$, whose depths are $\ge1$.
    Assume that $\Delta$ and $\sha$ are compatible for all words $\fa,\fb \in \langle \Gamma \rangle$ with $\depth(\fa)+\depth(\fb)<\depth(\fu)+ \depth(\fv)$. Further, we suppose that
    $\Delta(x_u\fu_- \sha x_v \fv_-) =\Delta(x_u\fu_-)\sha \Delta(x_v \fv_-)$.

Then we have \[ \Delta(\fu \sha \fv) = \Delta(\fu) \sha \Delta(\fv). \]
   \end{lemma}

    \begin{proof}
        Let \begin{align*}
            &\Delta (x_u )  = 1 \otimes x_u + \sum a_u \otimes b_u,&
            &\Delta (x_v )  = 1 \otimes x_v + \sum a_v \otimes b_v, \\
            &\Delta( \fu_-)  =  \sum a_{\fu_-} \otimes b_{\fu_-},&
            &\Delta( \fv_-)  =  \sum a_{\fv_-} \otimes b_{\fv_-},
        \end{align*}
        with $x_u, x_v, a_u, a_v \in \Sigma$, $b_u, b_v \in \langle \Sigma\rangle$ and $a_{\fu_-} , b_{\fu_-}, a_{\fv_-}, b_{\fv_-} \in \langle \Gamma \rangle$. We have
        \begin{align*}
            \Delta (x_{u, \alpha} ) & = 1 \otimes x_{u, \alpha} + \sum a_{u, \alpha} \otimes b_u, \\
            \Delta (x_{v, \beta} ) & = 1 \otimes x_{v, \beta} + \sum a_{v, \beta} \otimes b_v,
        \end{align*}
        where $a_{u, \alpha} = \varphi_\alpha (a_u)$, $a_{v, \beta} = \varphi_\beta(a_v)$. Note that
        \begin{align*}
            \Delta(x_{u, \alpha} \fu_-) & = 1 \otimes x_{u, \alpha} \fu_- +  \sum a_{u, \alpha}  a_{\fu_-}\otimes (b_u \sha b_{\fu_-}), \\
            \Delta(x_{v, \beta} {\fv_-}) & = 1 \otimes x_{v, \beta} \fv_- +  \sum a_{v, \beta}  a_{\fv_-} \otimes (b_v \sha b_{\fv_-}).
        \end{align*} Therefore, we have
        \begin{align}\label{eqn: step 2, delta shuffle delta}
             & \Delta(\fu)\sha \Delta(\fv) = \Delta(x_{u, \alpha} \fu_-) \shuffle \Delta( x_{v, \beta} \fv_-) \\
             =&
            \left( 1 \otimes x_{u, \alpha} \fu_-
            + \sum  a_{u, \alpha} a_{\fu_-} \otimes (b_u \shuffle b_{\fu_-})
            \right)  \shuffle \left(
            1 \otimes x_{v, \beta} \fv_- +  \sum a_{v, \beta}  a_{\fv_-} \otimes (b_v \sha b_{\fv_-})
            \right) \notag \\
             =&
            1 \otimes (x_{u, \alpha}\fu_- \sha x_{v, \beta}\fv_-)
            +
            \sum a_{v, \beta}  a_{\fv_-} \otimes (x_{u, \alpha} \fu_- \sha b_v \sha b_{\fv_-}) \notag \\
             &+
            \sum  a_{u, \alpha} a_{\fu_-} \otimes (x_{v, \beta}\fv_- \sha b_u \shuffle b_{\fu_-}) 
            +
            \sum  (a_{u, \alpha} a_{\fu_-} \sha a_{v, \beta} a_{\fv_-}) \otimes (b_u \shuffle b_v \sha b_{\fu_-}\sha b_{\fv_-}).
            \notag
        \end{align}
        On the other hand, from Lemma \ref{lem: distributive} and Lemma~\ref{lem: triangleright formulas}
        \begin{align*}
            \fu \sha \fv & =
            \fu \triangleright \fv
            + \fv \triangleright \fu
            + \varphi_{\alpha\beta}(x_u\fu_- \diamond x_v \fv_-),
        \end{align*}
        thus we have
        \begin{align}\label{eqn: step 2, delta of shuffle}
            & \Delta( \fu \sha \fv) = S_1 + S_2 + S_3
        \end{align}
        where
        \begin{align*}
            S_1 &= \Delta( \fu \triangleright \fv  )
            = \Delta( x_{u, \alpha}(\fu_- \sha \fv)), \\
            S_2 &= \Delta(\fv \triangleright \fu )
            = \Delta( x_{v, \beta} ( \fv_- \sha \fu )), \\
            S_3 & = \Delta(\varphi_{\alpha\beta}(x_u\fu_- \diamond x_v \fv_-)) \\
            & = (\varphi_{\alpha\beta}\otimes \Id) \Delta(x_u \fu_- \diamond x_v \fv_-)+  ( \Id  \otimes \varphi_{\alpha\beta} - \varphi_{\alpha\beta} \otimes\Id ) (1\otimes (x_u \fu_- \diamond x_v \fv_-)).
        \end{align*}

        We calculate $S_1$ and $S_2$. Since $\depth(\fu_-)+ \depth(\fv)<\depth(\fu)+ \depth(\fv)$, the induction hypothesis yields
        \begin{align*}
            \Delta(\fu_- \sha \fv) & = \Delta(\fu_-) \sha \Delta(x_{v, \beta}\fv_-) \\
            & = \left(\sum a_{\fu_-} \otimes b_{\fu_-}\right)
            \sha
            \left(
            1 \otimes x_{v, \beta}\fv_- + \sum a_{v, \beta} a_{\fv_-} \otimes (b_v \sha b_{\fv_-})
            \right)\\
            & = \sum a_{\fu_-} \otimes (b_{\fu_-} \sha x_{v, \beta} \fv_-) +
            \sum (a_{\fu_-} \sha a_{v, \beta} a_{\fv_-} )\otimes (b_v \sha b_{\fu_-} \sha b_{\fv_-}),
        \end{align*}
        so
        \begin{align*}
            S_1 & = 1 \otimes (x_{u, \alpha} (\fu_- \sha \fv)) +  \sum a_{u, \alpha} a_{\fu_-} \otimes (b_u \sha b_{\fu_-} \sha x_{v, \beta} \fv_-) \\
            & \quad + \sum a_{u, \alpha}(a_{\fu_-} \sha a_{v, \beta} a_{\fv_-} )\otimes (b_u \sha  b_v \sha b_{\fu_-} \sha b_{\fv_-}).
        \end{align*}

        A similar calculation yields
        \begin{align*}
            S_2 & = 1 \otimes (x_{v, \beta} (\fu \sha \fv_-)) +  \sum a_{v, \beta} a_{\fv_-} \otimes (b_v \sha b_{\fv_-} \sha x_{u, \alpha} \fu_-) \\
            & \quad + \sum a_{v, \beta}(a_{\fv_-} \sha a_{u, \alpha} a_{\fu_-} )\otimes (b_u \sha  b_v \sha b_{\fu_-} \sha b_{\fv_-}).
        \end{align*}

        What we want to show is that two equations~\eqref{eqn: step 2, delta shuffle delta} and \eqref{eqn: step 2, delta of shuffle} coincide, that is, $X:=\Delta(  \fu \sha \fv) -  \Delta(\fu) \sha \Delta(\fv ) = 0$.

        First, we collect the terms of the form $1\otimes \square$. The equation for the right tensorands of such terms in $X$ is
        \begin{align*}
            & x_{u, \alpha}( \fu_- \sha \fv) + x_{v, \beta}(\fu \sha\fv_-) + \varphi_{\alpha \beta } (x_u \fu_- \diamond x_v \fv_-) -  (x_u \fu_- \diamond x_v \fv_-) - \fu \sha\fv \\
            & = x_{u, \alpha}( \fu_- \sha \fv) + x_{v, \beta}(\fu \sha \fv_-) + \varphi_{\alpha \beta } (x_u \fu_- \diamond x_v \fv_-) -  (x_u \fu_- \diamond x_v \fv_-) \\
            & \quad - \fu \triangleright \fv
            - \fv \triangleright \fu
            - \fu \diamond \fv\\
            &= -  (x_u \fu_- \diamond x_v \fv_-)
        \end{align*}
        by Lemma~\ref{lem: distributive}, that is, the terms of the form $1\otimes \square$ in $X$ are collected as
        \[ - 1\otimes (x_u \fu_- \diamond x_v \fv_-).\]

        Next we collect the terms of the form $\sum \square \otimes (b_u \sha b_v \sha b_{\fu_-} \sha b_{\fv_-})$. For fixed indices, the equation for the left tensorands of such terms in $X$ is
        \begin{align*}
        & a_{u, \alpha} (a_{\fu_-} \sha a_{v, \beta}a_{\fv_-}) + a_{v, \beta}(a_{\fv_-} \sha a_{u, \alpha}a_{\fu_-}) - (a_{u, \alpha}a_{\fu_-} \sha a_{v, \beta}a_{\fv_-}) \\
        & = - \varphi_{\alpha \beta}(a_u a_{\fu_-} \diamond a_v a_{\fv_-}),
        \end{align*}
        by expanding the shuffle product of the third term and reducing it. That is, the terms of such a form in $X$ are collected as
        \[ -(\varphi_{\alpha \beta}\otimes \Id) \sum (a_u a_{\fu_-} \diamond a_v a_{\fv_-})\otimes (b_u \sha b_v \sha b_{\fu_-} \sha b_{\fv_-}). \]
        Some of the remaining terms cancel each other. One can check that
        \begin{align*}
            X =& ( \varphi_{\alpha\beta} \otimes \Id)\Delta( x_u \fu_- \diamond x_v \fv_-)  - 1\otimes (x_u \fu_- \diamond x_v \fv_-) \\
            & -(\varphi_{\alpha \beta}\otimes \Id) \sum (a_u a_{\fu_-} \diamond a_v a_{\fv_-})\otimes (b_u \sha b_v \sha b_{\fu_-} \sha b_{\fv_-}),
        \end{align*}
        which vanishes when
        \begin{equation}\Delta( x_u\fu_- \diamond x_v \fv_-)  = 1\otimes (x_u \fu_- \diamond x_v \fv_-) +\sum (a_u a_{\fu_-} \diamond a_v a_{\fv_-})\otimes (b_u \sha b_v \sha b_{\fu_-} \sha b_{\fv_-}).\label{eqn: compatibility induction step eqn1}\end{equation}

        Recall the assumption of the lemma, $\Delta(x_u \fu_- \sha x_v \fv_-) = \Delta(x_u \fu_-)\sha \Delta(x_v \fv_-)$. Expanding this by using definitions of $\sha$, $\diamond$ and calculating, we can conclude the desired Equation \eqref{eqn: compatibility induction step eqn1}. Thus the lemma holds.
    \end{proof}

    \begin{proposition}\label{prop: compatibility main induction step}
        Let $x_u \fu$, $x_v \fv$ be two words in $\langle \Gamma \rangle$, whose first letters are with trivial character $1 \in \mathbb{F}_q^*$. If $\Delta$ and $\sha$ are compatible for all words $\fa,\fb \in \langle \Gamma \rangle$ with $\depth(\fa)+\depth(\fb)<\depth(x_u \fu)+ \depth(x_v \fv)$, then
        \[ \Delta(x_u \fu \sha x_v \fv) = \Delta(x_u \fu)\sha \Delta(x_v \fv). \]
    \end{proposition}
    The proof closely follows that of \cite[Proposition 7.4]{IKLNDP23}, with the only difference being that the induction step is now based on depth rather than weight induction. For the sake of clarity, we include a complete proof.

    Note that the proof is valid for $\fu=1$ or $\fv=1$ cases. In those cases, some sum presentation is vacuous, i.e., if $\fu=1$, the indexed set $\{(a_\fu, b_\fu)\}$ is the empty set in the sum presentation of $\Delta(\fu) = 1 \otimes \fu + \sum a_{\fu} \otimes b_\fu$.

\begin{proof}[Proof of Proposition \ref{prop: compatibility main induction step}]
    By the definition of the $\shuffle$ and Lemma~\ref{lem: triangleright formulas}, we have
\begin{align} \label{eq: expansion 2 MZV}
& x_u \fu \shuffle x_v \fv \\
&=(x_u \fu) \triangleright (x_v \fv)+(x_v \fv) \triangleright (x_u \fu)+(x_u \fu) \diamond (x_v \fv) \notag \\
&=x_v(x_u \fu\shuffle \fv) + x_u(\fu\shuffle (x_v \fv))+(x_u \diamond x_v) \triangleright (\fu\shuffle \fv) \notag \\
&=x_v(x_u \fu\shuffle \fv) + x_u(\fu\shuffle (x_v \fv))+x_{u+v}(\fu\shuffle \fv)+ \sum_{i+j=u+v} \Delta^j_{u,v} x_{i} (x_j\shuffle \fu\shuffle \fv). \notag
\end{align}
Therefore, we get
\begin{align} \label{eq:compatibility b MZV}
 \Delta&(x_u \fu \shuffle x_v \fv)-\Delta(x_u \fu)\shuffle \Delta(x_v \fv) \\
=& \Delta(x_v(x_u \fu\shuffle \fv)) + \Delta(x_u(\fu\shuffle x_v \fv))+\Delta(x_{u+v}(\fu\shuffle \fv)) \notag \\
&+ \sum_{i+j=u+v} \Delta^j_{u,v} \Delta(x_{i} (x_j\shuffle \fu\shuffle \fv)) -\Delta(x_u \fu)\shuffle \Delta(x_v \fv). \notag
\end{align}

We calculate each term on the right hand side to conclude that this vanishes. We put
\begin{align*}
\Delta(\fu) =1 \otimes \fu+\sum a_\fu \otimes b_\fu, \text{ and }
\Delta(\fv) =1 \otimes \fv+\sum a_\fv \otimes b_\fv.
\end{align*}
Also for all $j \in \bN$, we put
\begin{align*}
\Delta(x_j) &=1 \otimes x_j+\sum a_j \otimes b_j.
\end{align*}
In particular,
\begin{align*}
\Delta(x_u) =1 \otimes x_u+\sum a_u \otimes b_u, \text{ and }
\Delta(x_v) =1 \otimes x_v+\sum a_v \otimes b_v.
\end{align*}

\noindent {\bf The first term $\Delta(x_v(x_u \fu\shuffle \fv))$.} \par
From the definition of $\Delta$, \begin{align}\label{eq: x_u u MZV}\Delta(x_u \fu)=1 \otimes x_u \fu+\sum a_u \otimes (b_u\shuffle \fu)+\sum (a_u \triangleright a_\fu) \otimes (b_u\shuffle b_\fu).\end{align} As $\depth(x_u \fu)+ \depth(\fv)<\depth(x_u \fu)+ \depth(x_v \fv)$, by the induction hypothesis we obtain
\begin{align*}
\Delta(x_u \fu\shuffle \fv)=& \Delta(x_u \fu)\shuffle \Delta(\fv) \\
=& 1 \otimes ((x_u \fu)\shuffle \fv)+\sum a_u \otimes (b_u\shuffle \fu\shuffle \fv)+\sum a_\fv \otimes ((x_u \fu)\shuffle b_\fv) \\
&+\sum (a_u \triangleright a_\fu) \otimes (b_u\shuffle b_\fu\shuffle \fv)+\sum (a_u \shuffle  a_\fv) \otimes (b_u\shuffle \fu\shuffle b_\fv) \\
&+\sum ((a_u \triangleright a_\fu)\shuffle a_\fv) \otimes (b_u\shuffle b_\fu\shuffle b_\fv).
\end{align*}
Thus
\begin{align} \label{eq:term 1b MZV}
\Delta(x_v(x_u \fu\shuffle \fv))=&1 \otimes x_v(x_u \fu\shuffle x_v)+ \sum a_v \otimes (b_v\shuffle (x_u \fu)\shuffle \fv) \\
&+\sum (a_v \triangleright a_u) \otimes (b_u\shuffle b_v\shuffle \fu\shuffle \fv) \notag \\
&+\sum (a_v \triangleright (a_u \triangleright a_\fu)) \otimes (b_u\shuffle b_v\shuffle b_\fu\shuffle \fv) \notag \\
&+\sum (a_v \triangleright a_\fv) \otimes (b_v\shuffle (x_u \fu)\shuffle b_\fv) \notag \\
&+\sum (a_v \triangleright (a_u\shuffle a_\fv)) \otimes (b_u\shuffle b_v\shuffle \fu\shuffle b_\fv) \notag \\
&+\sum (a_v \triangleright ((a_u \triangleright a_\fu)\shuffle a_\fv)) \otimes (b_u\shuffle b_v\shuffle b_\fu\shuffle b_\fv). \notag
\end{align}

\noindent {\bf The second term $\Delta(x_u(\fu\shuffle (x_v \fv))$.} \par Similarly, we get
\begin{align} \label{eq:term 2b MZV}
\Delta(x_u(\fu\shuffle (x_v \fv))=&1 \otimes x_u(\fu\shuffle x_v \fv)+ \sum a_u \otimes (b_u\shuffle \fu\shuffle (x_v \fv)) \\
&+\sum ((a_u \triangleright a_v) \otimes (b_u\shuffle b_v\shuffle \fu\shuffle \fv) \notag \\
&+\sum (a_u \triangleright (a_v \triangleright a_\fv)) \otimes (b_u\shuffle b_v\shuffle \fu\shuffle b_\fv) \notag \\
&+\sum (a_u \triangleright a_\fu) \otimes (b_u\shuffle b_\fu\shuffle (x_v \fv)) \notag \\
&+\sum (a_u \triangleright (a_\fu\shuffle a_v)) \otimes (b_u\shuffle b_\fu\shuffle b_v\shuffle \fv) \notag \\
&+\sum (a_u \triangleright (a_\fu\shuffle (a_v \triangleright a_\fv)) \otimes (b_u\shuffle b_v\shuffle b_\fu\shuffle b_\fv). \notag
\end{align}

\noindent {\bf The third term $\Delta(x_{u+v}(\fu\shuffle \fv))$.} \par

We put
	\[ \Delta(\fu\shuffle \fv)=1 \otimes (\fu\shuffle \fv)+\sum a_{\fu\shuffle \fv} \otimes b_{\fu\shuffle \fv}. \]
As $\depth(\fu)+ \depth(\fv)<\depth(x_u \fu)+ \depth(x_v \fv)$, the induction hypothesis implies that
\begin{align*}
\Delta(\fu\shuffle \fv)=\Delta(\fu)\shuffle \Delta(\fv) = \left(1 \otimes \fu+\sum a_\fu \otimes b_\fu \right)\shuffle \left(1 \otimes \fv+\sum a_\fv \otimes b_\fv \right).
\end{align*}
Thus
\begin{align*}
1 \otimes (\fu\shuffle \fv)+\sum a_{\fu\shuffle \fv} \otimes b_{\fu\shuffle \fv} 
=& \left(1 \otimes \fu+\sum a_\fu \otimes b_\fu \right)\shuffle \left(1 \otimes \fv+\sum a_\fv \otimes b_\fv \right) \\
=& 1 \otimes (\fu\shuffle \fv)+\sum a_\fu \otimes (b_\fu\shuffle \fv) \\
&+\sum a_\fv \otimes (\fu\shuffle b_\fv)+\sum (a_\fu\shuffle a_\fv) \otimes (b_\fu\shuffle b_\fv),
\end{align*}
which implies
\begin{equation} \label{eq: fu fv MZV}
\sum a_{\fu\shuffle \fv} \otimes b_{\fu\shuffle \fv}=\sum a_\fu \otimes (b_\fu\shuffle \fv)+\sum a_\fv \otimes (\fu\shuffle b_\fv)+\sum (a_\fu\shuffle a_\fv) \otimes (b_\fu\shuffle b_\fv).
\end{equation}
Thus we have
\begin{align} \label{eq:term 3b MZV}
\Delta(x_{u+v}(\fu\shuffle \fv)) =&1 \otimes (x_{u+v}(\fu\shuffle \fv))+ \sum a_{u+v} \otimes (b_{u+v}\shuffle \fu\shuffle \fv) \\
&+\sum (a_{u+v} \triangleright a_{\fu\shuffle \fv}) \otimes (b_{u+v}\shuffle b_{\fu\shuffle \fv}). \notag
\end{align}

\noindent {\bf The fourth terms $\Delta(x_{i} (x_j\shuffle \fu\shuffle \fv))$ for all $0 \le i$, $j \le u+v$ with $i+j=u+v$.}  \par

As $\depth(\fu)+ \depth(x_j)<\depth(x_u \fu)+ \depth(x_v \fv)$, we have $\Delta(x_j \sha \fu) = \Delta(\fu \sha x_j) = \Delta(\fu)\sha \Delta(x_j)$ by induction hypothesis. Further by Lemma~\ref{lem: sha and depth}, $ x_j\sha \fu$ can be represented as the sum of words with depth less than or equal to $\depth(x_j \fu)$. Thus $\Delta(x_j \sha \fu \sha \fv) = \Delta(x_j \sha \fu) \sha \Delta(\fv) = \Delta(x_j)\sha\Delta(\fu)\sha\Delta(\fv)= \Delta(x_j)\sha\Delta(\fu\sha\fv)$ by the induction hypothesis again.

Therefore,
\begin{align*}
\Delta(x_j\shuffle \fu\shuffle \fv) 
=&\Delta(x_j)\shuffle \Delta(\fu\shuffle \fv) \\
=& \left(1 \otimes x_j+\sum a_j \otimes b_j \right)\shuffle \left(1 \otimes (\fu\shuffle \fv)+\sum a_{\fu\shuffle \fv} \otimes b_{\fu\shuffle \fv} \right) \\
=& 1 \otimes (x_j\shuffle \fu\shuffle \fv)+\sum a_j \otimes (b_j\shuffle \fu\shuffle \fv) \\
&+\sum a_{\fu\shuffle \fv} \otimes (x_j\shuffle b_{\fu\shuffle \fv})+\sum (a_j\shuffle a_{\fu\shuffle \fv}) \otimes (b_j\shuffle b_{\fu\shuffle \fv}).
\end{align*}
This leads to
\begin{align}\label{eq:term 4b MZV}
 \Delta(x_{i} (x_j\shuffle \fu\shuffle \fv)) 
=& 1 \otimes (x_{i}(x_j\shuffle \fu\shuffle \fv))+\sum a_{i} \otimes (b_{i}\shuffle x_j\shuffle \fu\shuffle \fv)  \\
&+\sum (a_{i} \triangleright a_j) \otimes (b_{ i}\shuffle b_j\shuffle \fu\shuffle \fv) \notag \\
&+\sum (a_{ i} \triangleright a_{\fu\shuffle \fv}) \otimes (b_{i }\shuffle x_j\shuffle b_{\fu\shuffle \fv}) \notag \\
&+\sum (a_{ i} \triangleright (a_j\shuffle a_{\fu\shuffle \fv})) \otimes (b_{ i}\shuffle b_j\shuffle b_{\fu\shuffle \fv}) \notag.
\end{align}

\noindent {\bf The last term $\Delta(x_u \fu)\shuffle \Delta(x_v \fv)$.} \par

Note that $\Delta(x_u \fu)$ is given by \eqref{eq: x_u u MZV}. A parallel calculation gives
\begin{equation*}
\Delta(x_v \fv)=1 \otimes x_v \fv+\sum a_v \otimes (b_v\shuffle \fv)+\sum (a_v \triangleright a_\fv) \otimes (b_v\shuffle b_\fv).
\end{equation*}
Thus
\begin{align} \label{eq:term 5b MZV}
& \Delta(x_u \fu)\shuffle \Delta(x_v \fv) \\
= &\left(1 \otimes x_u \fu+\sum a_u \otimes (b_u\shuffle \fu)+\sum (a_u \triangleright a_\fu) \otimes (b_u\shuffle b_\fu) \right) \notag \\
&\shuffle \left(1 \otimes x_v \fv+\sum a_v \otimes (b_v\shuffle \fv)+\sum (a_v \triangleright a_\fv) \otimes (b_v\shuffle b_\fv)\right) \notag \\
=& 1 \otimes (x_u \fu\shuffle x_v \fv) \notag \\
&+ \sum a_v \otimes ((x_u \fu)\shuffle b_v\shuffle \fv)+\sum (a_v \triangleright a_\fv) \otimes ((x_u \fu)\shuffle b_v\shuffle b_\fv) \notag \\
&+ \sum a_u \otimes (b_u\shuffle \fu\shuffle (x_v \fv))+\sum (a_u\shuffle a_v) \otimes (b_u\shuffle b_v\shuffle \fu\shuffle \fv) \notag \\
&+\sum (a_u\shuffle (a_v \triangleright a_\fv)) \otimes (b_u\shuffle \fu\shuffle b_v\shuffle b_\fv) \notag \\
&+ \sum (a_u \triangleright a_\fu) \otimes (b_u\shuffle b_\fu\shuffle (x_v \fv))+ \sum ((a_u \triangleright a_\fu)\shuffle a_v) \otimes (b_u\shuffle b_\fu\shuffle b_v\shuffle \fv) \notag \\
&+\sum ((a_u \triangleright a_\fu)\shuffle (a_v \triangleright a_\fv)) \otimes (b_u\shuffle b_\fu\shuffle b_v\shuffle b_\fv). \notag
\end{align}

Plugging the equations \eqref{eq:term 1b MZV}, \eqref{eq:term 2b MZV}, \eqref{eq:term 3b MZV}, \eqref{eq:term 4b MZV}, \eqref{eq:term 5b MZV} into \eqref{eq:compatibility b MZV} yields
\begin{align*}
\Delta(x_u \fu \shuffle x_v \fv)-\Delta(x_u \fu)\shuffle \Delta(x_v \fv) =S_0-S_1-S_2-S_\fu-S_\fv+S_3+S_4.
\end{align*}
Here the sums $S_0$, $S_1$, $S_2$, $S_3$, $S_4$, $S_\fu$ and $S_\fv$ are given as follows:
\begin{align*}
S_0=& 1\otimes x_v(x_u \fu\shuffle \fv) + 1 \otimes x_u(\fu\shuffle (x_v \fv))+1 \otimes x_{u+v}(\fu\shuffle \fv) \\
&+ \sum_{i+j=u+v} \Delta^j_{u,v} 1 \otimes x_{i} (x_j\shuffle \fu\shuffle \fv)-1 \otimes ((x_u \fu)\shuffle (x_v \fv)). \\
S_1=&\sum (a_u\shuffle a_v) \otimes (b_u\shuffle b_v\shuffle \fu\shuffle \fv)-\sum (a_u \triangleright a_v) \otimes (b_u\shuffle b_v\shuffle \fu\shuffle \fv) \\
&-\sum (a_v \triangleright a_u) \otimes (b_u\shuffle b_v\shuffle \fu\shuffle \fv). \\
S_2=&\sum ((a_u \triangleright a_\fu)\shuffle (a_v \triangleright a_\fv)) \otimes (b_u\shuffle b_\fu\shuffle b_v\shuffle b_\fv) \\
&-\sum (a_u \triangleright (a_\fu\shuffle (a_v \triangleright a_\fv))) \otimes (b_u\shuffle b_\fu\shuffle b_v\shuffle b_\fv) \\
&-\sum (a_v \triangleright ((a_u \triangleright a_\fu)\shuffle a_\fv) \otimes (b_u\shuffle b_\fu\shuffle b_v\shuffle b_\fv). \\
S_3=& \sum a_{u+v} \otimes (b_{u+v}\shuffle \fu\shuffle \fv) +\sum_{i+j=u+v} \Delta^j_{u,v} \sum a_{ i} \otimes (b_{ i}\shuffle x_j\shuffle \fu\shuffle \fv) \\
&+ \sum_{i+j=u+v} \Delta^j_{u,v} \sum (a_{i} \triangleright a_j) \otimes (b_{i}\shuffle b_j\shuffle \fu\shuffle \fv). \\
S_4=& \sum (a_{u+v} \triangleright a_{\fu\shuffle \fv}) \otimes (b_{u+v}\shuffle b_{\fu\shuffle \fv}) +\sum_{i+j=u+v} \Delta^j_{u,v} \sum (a_{i} \triangleright a_{\fu\shuffle \fv}) \otimes (b_{i}\shuffle x_j\shuffle b_{\fu\shuffle \fv}) \\
&+ \sum_{i+j=u+v} \Delta^j_{u,v} \sum (a_{i} \triangleright (a_j\shuffle a_{\fu\shuffle \fv})) \otimes (b_{i}\shuffle b_j\shuffle b_{\fu\shuffle \fv}),
\end{align*}
and
\begin{align*}
S_\fu=&\sum (a_u\shuffle (a_v \triangleright a_\fv)) \otimes (b_u\shuffle b_v\shuffle \fu\shuffle b_\fv)-\sum (a_u \triangleright (a_v \triangleright a_\fv)) \otimes (b_u\shuffle b_v\shuffle \fu\shuffle b_\fv) \\
&-\sum (a_v \triangleright (a_u\shuffle a_\fv)) \otimes (b_u\shuffle b_v\shuffle \fu\shuffle b_\fv), \\
S_\fv=&\sum ((a_u \triangleright a_\fu)\shuffle a_v) \otimes (b_u\shuffle b_v\shuffle b_\fu\shuffle \fv)-\sum (a_u \triangleright (a_\fu\shuffle  a_v)) \otimes (b_u\shuffle b_v\shuffle b_\fu\shuffle \fv) \\
&-\sum (a_v \triangleright (a_u \triangleright a_\fu)) \otimes (b_u\shuffle b_v\shuffle b_\fu\shuffle \fv).
\end{align*}

We claim that
\begin{enumerate}
\item $S_0=0$.
\item $S_1-S_3=0$.
\item $S_2+S_\fu+S_\fv-S_4=0$.
\end{enumerate}
If these claims are true, the proposition is also true.

We prove three claims below. The claim (1), $S_0=0$, is immediate from \eqref{eq: expansion 2 MZV}.

We show the claim (2). We will show that
\begin{align*}
S_1=S_3=\sum (a_u \diamond a_v)  \otimes (b_u\shuffle b_v\shuffle \fu\shuffle \fv).
\end{align*}
With Lemma~\ref{lem: triangleright formulas}, we have
\begin{align*}
S_1=&\sum (a_u\shuffle a_v) \otimes (b_u\shuffle b_v\shuffle \fu\shuffle \fv)-\sum (a_u \triangleright a_v) \otimes (b_u\shuffle b_v\shuffle \fu\shuffle \fv) \\
&-\sum (a_v \triangleright a_u) \otimes (b_u\shuffle b_v\shuffle \fu\shuffle \fv) \\
=& \sum (a_u\shuffle a_v- a_u \triangleright a_v-a_v \triangleright a_u) \otimes (b_u\shuffle b_v\shuffle \fu\shuffle \fv) \\
=& \sum (a_u \diamond a_v)  \otimes (b_u\shuffle b_v\shuffle \fu\shuffle \fv).
\end{align*}

Now we calculate $S_3$. Recall $x_u \diamond x_v = x_{u+v} + \sum_{i+j=u+v}\Delta_{u,v}^j x_i x_j$.
Lemma~\ref{lem: delta diamond} gives that
\[\Delta(x_u \diamond x_v) = 1\otimes ( x_u \diamond x_v) + \sum (a_u \diamond a_v) \otimes (b_u \sha b_v),\]
so
\begin{align*}
    &1\otimes x_{u+v} +\sum a_{u+v}\otimes b_{u+v} \\
    &
    +\sum_{i+j=u+v} \Delta_{u,v}^j \left(
        1\otimes x_i x_j + \sum a_i \otimes  (b_i\sha x_j)  + \sum (a_i \triangleright a_j) \otimes (b_i \sha b_j)
        \right)\\
        \notag
          =& 1\otimes ( x_u \diamond x_v) + \sum (a_u \diamond a_v) \otimes (b_u \sha b_v),\end{align*}
        therefore
        \begin{align}
        \label{eq: x_u x_v 2 MZV}
    & \sum a_{u+v}\otimes b_{u+v} + \sum_{i+j=u+v} \Delta_{u, v}^j \sum a_i\otimes (b_i \sha x_j) +
    \sum_{i+j=u+v} \Delta_{u, v}^j \sum (a_i\triangleright a_j) \otimes (b_i \sha b_j)\\
    \notag
    & = \sum (a_u \diamond a_v) \otimes (b_u \sha b_v).
\end{align}

We conclude that
\begin{align*}
S_3=& \sum a_{u+v} \otimes (b_{u+v}\shuffle \fu\shuffle \fv) \\
&+\sum_{i+j=u+v} \Delta^j_{u,v} \sum a_{i} \otimes (b_{i}\shuffle x_j\shuffle \fu\shuffle \fv) + \sum_{i+j=u+v} \Delta^j_{u,v} \sum (a_{i} \triangleright a_j) \otimes (b_{i}\shuffle b_j\shuffle \fu\shuffle \fv) \\
=& \sum (a_u \diamond a_v) \otimes (b_u\shuffle b_v\shuffle \fu\shuffle \fv)
\end{align*}
and thus $S_1=S_3$ as desired.

We now show the claim (3). Precisely, we show that
\begin{align*}
S_2+S_\fu+S_\fv=S_4=\sum ((a_u \diamond a_v) \triangleright a_{\fu\shuffle \fv}) \otimes (b_u\shuffle b_v\shuffle b_{\fu\shuffle \fv}).
\end{align*}
We have
\begin{align*}
S_2=&\sum ((a_u \triangleright a_\fu)\shuffle (a_v \triangleright a_\fv)) \otimes (b_u\shuffle b_\fu\shuffle b_v\shuffle b_\fv) \\
&-\sum (a_u \triangleright (a_\fu\shuffle (a_v \triangleright a_\fv))) \otimes (b_u\shuffle b_\fu\shuffle b_v\shuffle b_\fv) \\
&-\sum (a_v \triangleright ((a_u \triangleright a_\fu)\shuffle a_\fv) \otimes (b_u\shuffle b_\fu\shuffle b_v\shuffle b_\fv) \\
=&\sum ((a_u \triangleright a_\fu)\shuffle (a_v \triangleright a_\fv)-a_u \triangleright (a_\fu\shuffle (a_v \triangleright a_\fv))-a_v \triangleright ((a_u \triangleright a_\fu)\shuffle a_\fv)) \\
&\quad \quad \otimes (b_u\shuffle b_\fu\shuffle b_v\shuffle b_\fv) \\
=&\sum ((a_u \triangleright a_\fu) \diamond (a_v \triangleright a_\fv)) \otimes (b_u\shuffle b_\fu\shuffle b_v\shuffle b_\fv).
\end{align*}
Here the last equality is from Lemma~\ref{lem: triangleright formulas} and Proposition~\ref{prop: properties on triangle diamond sha} as
\begin{align*}
&(a_u \triangleright a_\fu)\shuffle (a_v \triangleright a_\fv) \\&=(a_u \triangleright a_\fu) \triangleright (a_v \triangleright a_\fv)+(a_v \triangleright a_\fv) \triangleright (a_u \triangleright a_\fu)+(a_u \triangleright a_\fu) \diamond (a_v \triangleright a_\fv) \\
&=a_u \triangleright (a_\fu\shuffle (a_v \triangleright a_\fv))+a_v \triangleright ((a_u \triangleright a_\fu)\shuffle a_\fv))+(a_u \triangleright a_\fu) \diamond (a_v \triangleright a_\fv).
\end{align*}

By Lemma \ref{lem: triangleright formulas} and Proposition \ref{prop: properties on triangle diamond sha} again,
\begin{align*}
S_\fu=&\sum (a_u\shuffle (a_v \triangleright a_\fv)) \otimes (b_u\shuffle b_v\shuffle \fu\shuffle b_\fv)-\sum (a_u \triangleright (a_v \triangleright a_\fv)) \otimes (b_u\shuffle b_v\shuffle \fu\shuffle b_\fv) \\
&-\sum (a_v \triangleright (a_u\shuffle a_\fv)) \otimes (b_u\shuffle b_v\shuffle \fu\shuffle b_\fv) \\
=&\sum (a_u\shuffle (a_v \triangleright a_\fv)-a_u \triangleright (a_v \triangleright a_\fv)-a_v \triangleright (a_u\shuffle a_\fv)) \otimes (b_u\shuffle b_v\shuffle \fu\shuffle b_\fv) \\
=&\sum (a_u \diamond (a_v \triangleright a_\fv)) \otimes (b_u\shuffle b_v\shuffle \fu\shuffle b_\fv),
\end{align*}
and
\begin{align*}
S_\fv=&\sum ((a_u \triangleright a_\fu)\shuffle a_v) \otimes (b_u\shuffle b_v\shuffle b_\fu\shuffle \fv)-\sum (a_u \triangleright (a_\fu\shuffle  a_v)) \otimes (b_u\shuffle b_v\shuffle b_\fu\shuffle \fv) \\
&-\sum (a_v \triangleright (a_u \triangleright a_\fu)) \otimes (b_u\shuffle b_v\shuffle b_\fu\shuffle \fv) \\
=&\sum ((a_u \triangleright a_\fu) \diamond a_v) \otimes (b_u\shuffle b_v\shuffle b_\fu\shuffle \fv).
\end{align*}

With the equations for $S_2$, $S_\fu$, $S_\fv$ and \eqref{eq: fu fv MZV}, again we can apply Lemma~\ref{lem: triangleright formulas} and Proposition~\ref{prop: properties on triangle diamond sha} to have
\begin{align*}
S_2+S_\fu+S_\fv=\sum ((a_u \diamond a_v) \triangleright a_{\fu\shuffle \fv}) \otimes (b_u\shuffle b_v\shuffle b_{\fu\shuffle \fv}).
\end{align*}

We now consider the term $S_4$. We have
\begin{align*}
S_4=& \sum (a_{u+v} \triangleright a_{\fu\shuffle \fv}) \otimes (b_{u+v}\shuffle b_{\fu\shuffle \fv}) +\sum_{i+j=u+v} \Delta^j_{u,v} \sum (a_{i} \triangleright a_{\fu\shuffle \fv}) \otimes (b_{i}\shuffle x_j\shuffle b_{\fu\shuffle \fv}) \\
&+ \sum_{i+j=u+v} \Delta^j_{u,v} \sum (a_{i} \triangleright (a_j\shuffle a_{\fu\shuffle \fv})) \otimes (b_{i}\shuffle b_j\shuffle b_{\fu\shuffle \fv}) \\
=& \sum (a_{u+v} \triangleright a_{\fu\shuffle \fv}) \otimes (b_{u+v}\shuffle b_{\fu\shuffle \fv}) +\sum_{i+j=u+v} \Delta^j_{u,v} \sum (a_{i} \triangleright a_{\fu\shuffle \fv}) \otimes (b_{i}\shuffle x_j\shuffle b_{\fu\shuffle \fv}) \\
&+ \sum_{i+j=u+v} \Delta^j_{u,v} \sum ((a_{i} \triangleright a_j) \triangleright a_{\fu\shuffle \fv}) \otimes (b_{i}\shuffle b_j\shuffle b_{\fu\shuffle \fv}) \\
=& \sum ((a_u \diamond a_v) \triangleright a_{\fu\shuffle \fv}) \otimes (b_u\shuffle b_v\shuffle b_{\fu\shuffle \fv}).
\end{align*}
Here the second equality follows from Proposition \ref{prop: properties on triangle diamond sha}. The third one is a direct consequence of \eqref{eq: x_u x_v 2 MZV}.
This proves $S_2+S_\fu+S_\fv=S_4$ as claimed.

Now Proposition~\ref{prop: compatibility main induction step} is proven with the claims (1), (2), and (3), as
\begin{align*}
\Delta(x_u \fu \shuffle x_v \fv)-\Delta(x_u \fu)\shuffle \Delta(x_v \fv)&=S_0-S_1-S_2-S_\fu-S_\fv+S_3+S_4 \\
&=S_0-(S_1-S_3)-(S_2+S_\fu+S_\fv-S_4) \\
&=0.
\end{align*}
\end{proof}

\begin{proof}[Proof of Theorem~\ref{thm: main statement on compatibility}] By combining Lemma~\ref{lemma: compatibility first character reduce} and Proposition~\ref{prop: compatibility main induction step} we conclude that, if Theorem~\ref{thm: main statement on compatibility} holds for all words $(\fa,\fb) \in \langle \Gamma \rangle \times \langle \Gamma \rangle$ with $\depth(\fa)+\depth(\fb)<n$ for some $n \in \mathbb N$ then the theorem holds for all words $(x_u \fu, x_v \fv)$ such that $\depth(x_u \fu)+\depth(x_v \fv)=n$, and thus for $(x_{u, \alpha} \fu, x_{v, \beta}\fv)$ such that $\depth(x_{u, \alpha} \fu)+\depth(x_{v, \beta} \fv)=n$, for all $u, v \in \mathbb{N}$ and $\alpha, \beta \in \mathbb{F}_q^*$. Since the theorem holds for the initial cases $(1, \fb)$ and $(\fa, 1)$ for all $\fa, \fb \in \langle \Gamma \rangle$ with the arbitrary depths, the theorem holds for all $\fu, \fv \in \langle \Gamma \rangle$.
\end{proof}

\subsection{Coassociativity}  \label{sec: coassociativity} ${}$\par

Next we want to prove:

\begin{theorem} \label{thm: coassociativity}
Let $\fu \in \langle \Gamma \rangle$. Then we have
	\[ (\Id \otimes \Delta) \Delta(\fu)=(\Delta \otimes \Id)\Delta(\fu). \]
\end{theorem}
\begin{proof}
The proof is formed into two parts as follows:
\begin{itemize}
\item Part (1): if the coassociativity holds for $\fu = x_{u} \fv$, then coassociativity holds for $\fu = x_{u, \varepsilon} \fv$, for $\fv \in\langle \Gamma\rangle $ with $\depth(\fv)\ge0$, and
\item Part (2): if the coassociativity holds for all $\fu$ with $\depth(\fu) \le n$, then the coassociativity holds for $\fu = x_{u} \fv$ for $u \in \mathbb{N}$ and $\depth \fv = n$.
\end{itemize}

Note that the initial case (i.e., when $\fu = x_{u, \varepsilon}$) is established by the coassociativity in~$\mathfrak{C}$ and Part~(1).

    For Part (1), we assume that $(\Id\otimes \Delta) \Delta(x_u \fv)  = (\Delta\otimes \Id)\Delta(x_u \fv)$. Let
    \begin{align*}
    \Delta(x_u) = 1\otimes x_u + \sum a_u \otimes b_u, \quad\quad
    \Delta(\fv)   = \sum a_\fv \otimes b_\fv.
    \end{align*}
Here $a_u \in \Sigma$, $b_u \in \langle \Sigma \rangle$, and $a_\fv, b_\fv \in \langle \Gamma \rangle$. Thus
	\[ \Delta( x_u \fv) = 1\otimes x_u \fv + \sum a_u a_\fv \otimes (b_u \sha b_\fv), \]
and
    \begin{align*}
        (\Id\otimes \Delta) \Delta( x_u \fv) & = 1\otimes 1\otimes x_u \fv + \sum 1\otimes a_u a_\fv \otimes (b_u \sha b_\fv)  + \sum a_u a_\fv \otimes \Delta(b_u \sha b_\fv), \\
        (\Delta \otimes \Id)\Delta(x_u \fv) & = 1\otimes1\otimes x_u \fv + \sum \Delta(a_u a_\fv) \otimes (b_u\sha b_\fv).
    \end{align*}
    Then the induction hypothesis that is the coassociativity for $x_u \fv$ is equivalent to
    \begin{align}
        \sum 1\otimes a_u a_\fv \otimes (b_u \sha b_\fv)  + \sum a_u a_\fv \otimes \Delta(b_u \sha b_\fv)= \sum \Delta(a_u a_\fv) \otimes (b_u\sha b_\fv). \label{eqn: coasso part 1 assumption}
    \end{align}

    Now let $\fu = x_{u, \varepsilon} \fv$. We have
    \begin{align*}
    \Delta(x_{u, \varepsilon}) & = 1\otimes x_{u, \varepsilon} + \sum a_{u, \varepsilon} \otimes b_u
    \end{align*}
    where $a_{u, \varepsilon} = \varphi_\varepsilon (a_u)$. By definition, we have
    \begin{align*}
        \Delta( \fu) = 1\otimes \fu + \sum a_{u, \varepsilon}a_\fv \otimes (b_u \sha b_\fv).
    \end{align*}
    One can calculate
    \begin{align*}
        (\Id \otimes \Delta) \Delta(\fu)  =& 1\otimes \Delta(\fu) + \sum a_{u, \varepsilon}a_\fv \otimes \Delta(b_u \sha b_\fv)\\
        &= 1\otimes 1\otimes \fu + \sum 1\otimes a_{u, \varepsilon}a_\fv \otimes (b_u \sha b_\fv) + \sum a_{u, \varepsilon}a_\fv \otimes \Delta(b_u \sha b_\fv),\\
        (\Delta \otimes \Id)\Delta(\fu)  =& 1\otimes 1\otimes \fu + \sum \Delta(a_{u, \varepsilon} a_\fv )\otimes (b_u \sha b_\fv) \\
         =& 1\otimes 1\otimes \fu +
        \sum
        (\varphi_\varepsilon \otimes \Id)\Delta(a_u a_\fv)
        \otimes (b_u \sha b_\fv) \\
         &+
        \sum
        1\otimes a_{u, \varepsilon} a_\fv
        \otimes (b_u \sha b_\fv) -
        \sum
        1\otimes a_u a_\fv
        \otimes (b_u \sha b_\fv).
    \end{align*}
    Thus the coassociativity holds for $x_{u, \varepsilon} \fv$ if
    \begin{align*}
        \sum
        (\varphi_\varepsilon \otimes \Id)\Delta(a_u a_\fv)
        \otimes (b_u \sha b_\fv)  =  \sum a_{u, \varepsilon}a_\fv \otimes \Delta(b_u \sha b_\fv) + \sum
        1\otimes a_u a_\fv
        \otimes (b_u \sha b_\fv),
    \end{align*}
    which is obtained by applying $\varphi_\varepsilon \otimes \Id \otimes \Id$ on Equation~\eqref{eqn: coasso part 1 assumption}. This proves that the coassociativity for $x_u \fv$ implies the coassociativity for $x_{u, \varepsilon} \fv$.

    We now prove Part (2). Let $n\ge1$. We assume that the coassociativity holds for all words with depth $\le n$. Let $\fu = x_{t, \varepsilon}\mathfrak{t}$ with depth $\le n$. We write
    \begin{align*}
        \Delta(x_{t, \varepsilon}) = 1\otimes x_{t, \varepsilon} +
        \sum a_{t, \varepsilon} \otimes b_t, \text{ and }
        \Delta(\mathfrak{t}) =  \sum a_\ft \otimes b_\mathfrak{t}.
    \end{align*}
where $a_{t, \varepsilon} \in \Gamma$ and $b_{t, \varepsilon}, a_\ft, b_\ft \in \langle \Gamma \rangle$. Then
	\[ \Delta(x_{t, \varepsilon} \mathfrak{t}) = 1\otimes x_{t, \varepsilon} + \sum a_{t, \varepsilon} a_{\mathfrak t} \otimes (b_t \sha b_{\mathfrak t}), \]
and the coassociativity for $x_{t, \varepsilon} \mathfrak{t}$ is reduced into
    \begin{align}\sum \left( \Delta(a_{t, \varepsilon} a_\mathfrak{t}) - 1\otimes a_{t, \varepsilon} a_\mathfrak{t} \right) \otimes (b_t \sha b_\mathfrak{t}) = \sum a_{t, \varepsilon} a_\mathfrak{t} \otimes (\Delta(b_t) \sha \Delta(b_\mathfrak{t})). \label{eqn: coasso part 2 assumption}
    \end{align}
    As a special case, let $\mathfrak{t} = 1$. Then $\Delta(\mathfrak {t}) = 1\otimes 1$, and
    \begin{align*}\sum \left( \Delta(a_{t, \varepsilon} ) - 1\otimes a_{t, \varepsilon}  \right) \otimes b_t = \sum a_{t, \varepsilon}\otimes \Delta(b_t).
    \end{align*}

    Now assume $\mathfrak{u} = x_u x_{v, \beta} \fv$ with depth $n+1$. Say
    \begin{align*}
    \Delta(x_u) &= 1\otimes x_u + \sum a_u \otimes b_u, \\
    \Delta(x_v) &= 1\otimes x_v + \sum a_v \otimes b_v, \quad \text{so } \Delta(x_{v, \beta}) = 1\otimes x_{v, \beta}+ \sum a_{u,\beta} \otimes b_v,\\
    \Delta(\fv) &  = \sum a_\fv \otimes b_\fv,
    \end{align*}
    where $a_u, a_v \in \Sigma$, $b_u, b_v \in \langle \Sigma \rangle$, $a_\fv, b_\fv \in \langle \Gamma \rangle$, and $a_{u, \beta} = \varphi_\beta (a_u)$. Then,
    \begin{align*}
        \Delta( x_{v, \beta} \fv) = 1\otimes x_{v, \beta}\fv + \sum a_{v, \beta}a_\fv \otimes (b_v \sha b_\fv).
    \end{align*}
It follows that
    \begin{align*}
        \Delta(x_u x_{v, \beta} \fv) = 1\otimes x_u  x_{v, \beta} \fv +
        \sum a_u\otimes (b_u \sha x_{v, \beta} \fv) +
        \sum a_u a_{v, \beta}a_\fv \otimes (b_u \sha b_v \sha b_\fv),
    \end{align*}
    thus
    \begin{align*}
        & (\Id\otimes \Delta)\Delta(x_u x_{v, \beta} \fv) \\
          =& 1\otimes 1\otimes x_u  x_{v, \beta} \fv +
        \sum 1\otimes  a_u\otimes (b_u \sha x_{v, \beta} \fv) +
        \sum 1\otimes  a_u a_{v, \beta}a_\fv \otimes (b_u \sha b_v \sha b_\fv)\\
        & + \sum a_u\otimes\Delta( b_u \sha x_{v, \beta} \fv) +
        \sum a_u a_{v, \beta}a_\fv \otimes \Delta (b_u \sha b_v \sha b_\fv),\\
        & (\Delta\otimes \Id)\Delta(x_u x_{v, \beta} \fv) \\
        =& 1\otimes 1 \otimes x_u x_{v, \beta} \fv +
        \sum \Delta(a_u) \otimes (b_u \sha x_{v, \beta} \fv) +
        \sum \Delta (a_u a_{v, \beta} a_\fv) \otimes (b_u \sha b_v \sha b_\fv),
    \end{align*}
    so what we want to show is that
    \begin{align} \label{eqn: coasso part2 mid-1}
        & \sum 1\otimes  a_u\otimes (b_u \sha x_{v, \beta} \fv) +
        \sum 1\otimes  a_u a_{v, \beta}a_\fv \otimes (b_u \sha b_v \sha b_\fv)\\
        &  + \sum a_u\otimes\Delta( b_u \sha x_{v, \beta} \fv) +
        \sum a_u a_{v, \beta}a_\fv \otimes \Delta (b_u \sha b_v \sha b_\fv), \notag\\
         =&
        \sum \Delta(a_u) \otimes (b_u \sha x_{v, \beta} \fv) +
        \sum \Delta (a_u a_{v, \beta} a_\fv) \otimes (b_u \sha b_v \sha b_\fv). \notag
    \end{align}

    Note that, by Equation~\eqref{eqn: coasso part 2 assumption}  we have
    \begin{align*}
        & \sum 1\otimes  a_u\otimes (b_u \sha x_{v, \beta} \fv) + \sum a_u\otimes\Delta( b_u \sha x_{v, \beta} \fv) -  \sum \Delta(a_u) \otimes (b_u \sha x_{v, \beta} \fv) \\
        & = \sum a_u\otimes\Delta( b_u \sha x_{v, \beta} \fv) - \left( \sum ( \Delta(a_u) - 1\otimes  a_u) \otimes b_u  \right) \sha  (1\otimes 1\otimes x_{v, \beta} \fv)
         \\
        & =  \sum a_u\otimes\Delta( b_u \sha x_{v, \beta} \fv) - \left( \sum a_u \otimes \Delta (b_u) \right)\sha  (1\otimes 1\otimes x_{v, \beta} \fv)\\
        & = \sum a_u\otimes\left( \Delta( b_u) \sha \Delta(x_{v, \beta} \fv)\right) -  \sum a_u \otimes \left( \Delta (b_u) \sha(1\otimes x_{v, \beta} \fv)\right)\\
        & =
        \sum  a_u \otimes \left( \Delta(b_u) \sha (\Delta(x_{v, \beta} \fv)-  1\otimes x_{v, \beta} \fv)
        \right)  \\
        & =\sum  a_u \otimes \left( \Delta(b_u) \sha (
        a_{v, \beta}a_\fv \otimes (b_v\sha b_\fv)
        )
        \right)
    \end{align*}
    so Equation~\eqref{eqn: coasso part2 mid-1} is equivalent to
    \begin{align} \label{eqn: coasso part2 mid-1-1}
        & \sum \Delta (a_u a_{v, \beta} a_\fv) \otimes (b_u \sha b_v \sha b_\fv) \\
        &=  \sum 1\otimes  a_u a_{v, \beta}a_\fv \otimes (b_u \sha b_v \sha b_\fv) + \sum a_u a_{v, \beta}a_\fv \otimes \Delta (b_u \sha b_v \sha b_\fv) \notag\\
        & \quad + \sum  a_u \otimes \left( \Delta(b_u) \sha (
        a_{v, \beta}a_\fv \otimes (b_v\sha b_\fv)
        )
        \right).\notag
    \end{align}

From now on, tracing the indices is important. When needed, we use the notation $\sum_{u}$ which stands for the sum over the pairs $(a_u, b_u)$, and similarly $\sum_{(v, \beta)}$, $\sum_\fv$, etc.

    We now fix $a_u$, $b_u$. The last two terms of \eqref{eqn: coasso part2 mid-1-1} become
    \begin{align*}
        & \sum a_u a_{v, \beta}a_\fv \otimes \left( \Delta (b_u) \sha \Delta(b_v) \sha \Delta(b_\fv)\right)  + \sum  a_u \otimes \left( \Delta(b_u) \sha (
        a_{v, \beta}a_\fv \otimes (b_v\sha b_\fv)
        )
        \right) \\
        & = (1\otimes \Delta(b_u)) \sha \left( \sum a_u a_{v, \beta}a_\fv \otimes \left( \Delta(b_v) \sha \Delta(b_\fv)\right)  + \sum  a_u \otimes
        a_{v, \beta}a_\fv \otimes (b_v\sha b_\fv)
        \right) \\
        & = (a_u \otimes 1 \otimes 1 ) \cdot \left( (1\otimes \Delta(b_u)) \sha \left( \sum a_{v, \beta}a_\fv \otimes \left( \Delta(b_v) \sha \Delta(b_\fv)\right)  + \sum  1 \otimes
        a_{v, \beta}a_\fv \otimes (b_v\sha b_\fv)
        \right)
        \right) \\
        & =(a_u \otimes 1 \otimes 1 ) \cdot \left( (1\otimes \Delta(b_u)) \sha \left(
        \sum \Delta(a_{v, \beta} a_{\fv}) \otimes (b_{v}\sha b_\fv)
        \right)
        \right).
    \end{align*}
    Recall that $\cdot$ is the concatenation. The last equality is from the induction hypothesis, Equation~\eqref{eqn: coasso part 2 assumption}. Equation~\eqref{eqn: coasso part2 mid-1-1} becomes
    \begin{align}\label{eqn: coasso part2 mid-2}
        & \sum \Delta (a_u a_{v, \beta} a_\fv) \otimes (b_u \sha b_v \sha b_\fv) \\
        &=  \sum 1\otimes  a_u a_{v, \beta}a_\fv \otimes (b_u \sha b_v \sha b_\fv)
        +\sum (a_u \otimes 1 \otimes 1 ) \cdot \left( (1\otimes \Delta(b_u)) \sha \left(
        \sum \Delta(a_{v, \beta} a_{\fv}) \otimes (b_{v}\sha b_\fv)
        \right)
        \right).
        \notag
    \end{align}

    Next we fix $a_u$, $a_{v, \beta}$, $a_\fv$. So the corresponding $b$'s are also fixed.
    Let
    \begin{align*}
    \Delta (a_u) = 1\otimes a_u + \sum_{(1)} c_{(1)} \otimes d_{(1)}, \quad
    \Delta (a_{v, \beta} a_\fv)  = \sum_{(2)} c_{(2)} \otimes d_{(2)}.
    \end{align*}
  Then, we have
    \begin{align*}
        \Delta(a_u a_{v, \beta} a_\fv) = 1\otimes a_u a_{v, \beta} a_\fv + \sum_{(1), (2)} c_{(1)}c_{(2)}\otimes (d_{(1)}\sha d_{(2)}).
    \end{align*}
The summand of $\sum \Delta (a_u a_{v, \beta} a_\fv) \otimes (b_u \sha b_v \sha b_\fv)$ for fixed $a_u$, $a_{v,\beta}$, $a_\fv$ is
    \begin{align*}
        \left( 1\otimes a_u a_{v, \beta} a_\fv + \left(\sum_{(1), (2)} c_{(1)}c_{(2)}\otimes (d_{(1)}\sha d_{(2)})\right) \right) \otimes (b_u \sha b_v \sha b_\fv).
    \end{align*}
Plugging this and the sum representation of $\Delta(a_{v, \beta} a_\fv)$ into the equation~\eqref{eqn: coasso part2 mid-2} yields
    \begin{align*}
        & \sum_{u, (v, \beta), \fv}
        \left(
        \sum_{(1), (2)} c_{(1)}c_{(2)}\otimes (d_{(1)}\sha d_{(2)})
        \right) \otimes (b_u \sha b_v \sha b_\fv)\\
        -&  \sum  a_u  c_{(2)} \otimes (\Delta(b_u) \sha (d_{(2)}  \otimes b_{v}\sha b_\fv)).
    \end{align*}

Finally, we fix $a_u$, $a_{v, \beta}$, $a_\fv$, $c_{(2)}$. The corresponding summands are
    \begin{align*}
        & \sum_{(1)} c_{(1)}c_{(2)}\otimes (d_{(1)}\sha d_{(2)}) \otimes (b_u \sha b_v \sha b_\fv) -  a_u c_{(2)} \otimes (\Delta(b_u) \sha (d_{(2)}  \otimes b_{v}\sha b_\fv)).
    \end{align*}
    By direct calculations, we get
    \begin{align*}
        & \sum_{(1)} c_{(1)}c_{(2)}\otimes (d_{(1)}\sha d_{(2)}) \otimes (b_u \sha b_v \sha b_\fv) -  a_u c_{(2)} \otimes (\Delta(b_u) \sha (d_{(2)}  \otimes b_{v}\sha b_\fv)) \\
        & =\left( \left( \sum_{(1)} c_{(1)}\otimes d_{(1)} \otimes b_u  -  a_u  \otimes \Delta(b_u)
        \right)
        \sha( 1\otimes d_{(2)} \otimes b_{v}\sha b_{\fv})\right)
        \cdot (c_{(2)} \otimes 1\otimes 1),
    \end{align*}
    whose inner sum is reduced to
    \begin{align*}
        & \sum_{(1)} c_{(1)}\otimes d_{(1)} \otimes b_u  -  a_u  \otimes 1  \otimes \Delta(b_u)  = (
        \Delta (a_u) - 1\otimes a_u
        )  \otimes b_u  -  a_u  \otimes 1  \otimes \Delta(b_u).
    \end{align*}
    Now we can apply Equation~\eqref{eqn: coasso part 2 assumption} for the sum of the summand over $a_u$ to have
    \begin{align*}
        \sum_u
        \left(\Delta (a_u)\otimes b_u - 1\otimes a_u \otimes b_u-  a_u  \otimes 1  \otimes \Delta(b_u)
        \right)  =0,
    \end{align*}
    that is, Equation~\eqref{eqn: coasso part2 mid-2} is true. Thus the induction step is proved.

    From the coassociativity of $\mathfrak{C}$, we have the coassociativity for $\fu \in \Sigma$, and then for $\fu \in \Gamma$ by Part (1). By applying Parts (2) and (1) repeatedly, we can cover arbitrary depth of $\fu \in \langle \Gamma \rangle$. Hence this completes the proof.
\end{proof}

\subsection{Hopf algebra structure on $\mathfrak{D}$} \ppar

Let the counit $\epsilon \colon \mathfrak{D} \to \Fq$ defined as $\epsilon(1) = 1$; $\epsilon(\mathfrak{u}) = 0$ otherwise.

We already checked that $(\mathfrak{D}, \shuffle, u, \Delta, \epsilon)$ is a bialgebra by Theorem~\ref{thm: commutative algebra}, Theorem~\ref{thm: main statement on compatibility}, and Theorem~\ref{thm: coassociativity}. By induction it is easily verified that $\sha$ and $\Delta$ preserve the weight, and the weight zero subset is isomorphic to $\mathbb{F}_q$. Therefore,  $(\mathfrak{D}, \shuffle, u, \Delta, \epsilon)$ is a connected graded bialgebra. By Theorem~\ref{prop: graded Hopf algebras}, we have the following theorem.
\begin{theorem} \label{thm: AMZV Hopf algebra}
$(\mathfrak{D}, \shuffle, u, \Delta, \epsilon)$  is a connected graded Hopf algebra of finite type over~$\mathbb{F}_q$.
\end{theorem}


\end{document}